# Weak backward error analysis for overdamped Langevin processes

Marie Kopec

*ENS Cachan Bretagne, campus de Ker Lann, avenue Robert Schumann, 35170 Bruz, France*

**Abstract**

We consider numerical approximations of overdamped Langevin stochastic differential equations by implicit methods. We show a weak backward error analysis result in the sense that the generator associated with the numerical solution coincides with the solution of a modified Kolmogorov equation up to high order terms with respect to the stepsize. This implies that every measure of the numerical scheme is close to a modified invariant measure obtained by asymptotic expansion. Moreover, we prove that, up to negligible terms, the dynamic associated with the implicit scheme considered is exponentially mixing.

## 1 Introduction

In [4], the authors give a weak backward error analysis for SDEs defined on the d-dimensional torus. The aim of this article is to extend the result of [4] to the overdamped Langevin process on $\mathbb{R}^d$.

In the last decades, backward error analysis has become a powerful tool to analyze the long time behavior of numerical schemes applied to evolution equations (see [10, 13, 20]). The main idea can be described as follows: Let us consider an ordinary differential equation of the form

$$\dot{y}(t) = f(y(t)),$$

where $f : \mathbb{R}^d \to \mathbb{R}^d$ is a smooth vector field, and denote by $\phi_t^f(y)$ the associated flow. By definition, a numerical method defines for a small time step $\delta$ an approximation $\Phi_\delta$ of the exact flow $\phi_\delta^f$: We have for bounded $y \in \mathbb{R}^d$, $\Phi_\delta(y) = \phi_\delta^f(y) + \mathcal{O}(\delta^{r+1})$ where $r$ is the order of the method.

The idea of backward error analysis is to show that $\Phi_\delta$ can be interpreted as the exact flow $\phi_\delta^{f_\delta}$ of a modified vector field defined as a series in powers of $\delta$

$$f_\delta = f + \delta^r f_r + \delta^{r+1} f_{r+1} + ...,$$

where $f_l$, $l \geq r$ are vector fields depending on the numerical method. In general, the series defining $f_\delta$ does not converge, but it can be shown that for bounded $y$, we have for arbitrary $N$

$$\Phi_\delta(y) = \phi_\delta^{f_\delta^N}(y) + C_N \delta^{N+1},$$

where $f_\delta^N$ is the truncated series:

$$f_\delta^N = f + \delta^r f_r + ... + \delta^N f_N.$$

Under some analytic assumptions, the constant $C_N \delta^{N+1}$ can be optimized in $N$, so that the error term in the previous equation can be made exponentially small with respect to $\delta$.

Such a result is very important and has many applications in the case where $f$ has some strong geometric properties, such as Hamiltonian (see [9, 10, 13, 18, 20]). In this situation, and under some compatibility conditions on the numerical method $\Phi_\delta$, the modified vector field $f_\delta$ inherits the structure of $f$.

More recently, these ideas have been extended in some situations to Hamiltonian PDEs : First in the linear case [3] and then in the semi linear case (nonlinear Schrödingier or wave equations), see [6, 7].

We want to use this approach for the stochastic differential equation:

$$dX(t) = -DV(X(t))dt + dW(t), \tag{1.1}$$

where $V$ and all its derivatives have polynomial growth. This equation describes the evolution of particles in a potential such that the corresponding canonical measure has several regions of high probability separated by low probability regions (see [14]). Moreover, the invariant measure associated to (1.1) represents the position of the particles at the equilibrium. It is interesting to have informations on its approximation. This equation can also describe the blood clotting dynamics [12].

In this work, we investigate the weak error which concerns the law of the solution. Let us recall that given a SDE in $\mathbb{R}^d$ of the form

$$dX = f(X)dt + g(X)dW, \tag{1.2}$$

discretized by an explicit Euler scheme $(X_p)$ with time step $\delta$, then, under assumptions on $f$, $g$ and $\phi : \mathbb{R}^d \to \mathbb{R}$ (see [16, 17, 24, 25]), the explicit Euler scheme $(X_p)$ has weak error order 1:

$$|\mathbb{E}(\phi(X_p)) - \mathbb{E}(\phi(X(p\delta)))| \leq c(\phi, T)\delta, \quad p = 0, ..., \lfloor T/\delta \rfloor, \quad T > 0.$$

Error estimates on long times for elliptic and hypoelliptic SDEs have already been proved, especially in the case of explicit scheme. In [23, 24], it is shown that for a sufficiently small time step, the explicit Euler scheme defines an ergodic process and that the invariant measure of the Euler scheme is close to the invariant measure of the SDE. In [26], under the assumption of the existence of a unique invariant measure associated to the SDE, Talay and Tubaro have shown that the weak error and the invariant measure associated to the Euler scheme can be expanded in powers of the time step $\delta$. The assumptions on $f$ and $g$ used in [23, 24, 26] are restrictive. Moreover, the results describe in these papers are only for explicit schemes.

In [11, 25], Higham, Mattingly, Stuart and Talay work with implicit schemes, but they need that the time step is small enough. In [11], under some assumptions, it is shown, in particular in the case of the overdamped Langevin equation, that for sufficiently small time step, two kinds of implicit schemes are ergodic processes. They also show that the invariant measures associated with theses schemes converge to the invariant measure of the overdamped Langevin equation. In [25], Talay studies stochastic Hamiltonian system. He shows the exponential convergence of the solution associated to the Kolmogorov equation and, for a sufficiently small time step $\delta$, an expansion with respect to $\delta$ of the invariant measure to the implicit scheme which is close to the invariant measure of the SDE.

In [16], a larger class of schemes is studied. It is shown that given an elliptic or hypoelliptic SDE defined on the d-dimensional torus, the ergodic averages provided by a class of implicit and explicit schemes are asymptotically close to the average of the invariant measure of the SDE. The authors also show an expansion in expectation of the invariant measure for any time-step.

In this paper, we work on $\mathbb{R}^d$. Moreover, in our case, $V$ and all its derivates are not bounded but have polynomial growth. The aim of this paper is, under assumptions on $V$, to show a weak backward error analysis result : We show an expansion with respect to the time step $\delta$ of $\mathbb{E}\phi(X_p)$ where $X_p$ is an implicit scheme. The idea to extend the backward error analysis to SDE has already be studied in [1, 4, 22]. In [1], the authors use this approach to construct new methods of weak order two to approximate stochastic differential equations. In [4], the authors study a SDE defined on the d-dimensional torus and its approximation by the explicit Euler scheme. They show, without restriction on the time step, an expansion of $\mathbb{E}\phi(X_p)$ where $X_p$ is the explicit Euler scheme. In [22], Shardlow consider SDE with additive noise ($g$ does not depend of $X$). He has shown that it is possible to build a modified SDE associated with the Euler scheme, but only at the first step, i.e. for $N = 2$. In this case, he is able to write down a modified SDE:

$$d\tilde{X} = \tilde{f}(\tilde{X})dt + \tilde{g}(\tilde{X})dW,$$

such that

$$\left|\mathbb{E}\big(\phi(X_p)\big) - \mathbb{E}\big(\phi(\tilde{X}(p\delta))\big)\right| \leq c_1(\phi, T)\delta^2, \ p = 0, ..., \lfloor T/\delta \rfloor, \ T > 0.$$

In this paper, we take the approach describes in [4]. We show that the generator associated with the process solution of the SDE coincides with the solution of a modified Kolmogorov equation up to high order terms with respect to the stepsize. It is known that given $\phi : \mathbb{R}^d \to \mathbb{R}$ and denoting by $X_x(t)$ the solution of the SDE (1.1) satisfying $X(0) = x$, the function $u$ defined for $t \geq 0$ and $x \in \mathbb{R}^d$ by $u(t, x) = \mathbb{E}(\phi(X_x(t)))$ satisfies the Kolmogorov equation

$$\partial_t u = Lu,$$

where $L$ is the Kolmogorov operator associated with the SDE.
We show that with the numerical solution, we can associate a modified Kolmogorov operator of the form

$$\mathcal{L}(\delta, x, \partial_x) = L(x, \partial_x) + \delta L_1(x, \partial_x) + \delta^2 L_2(x, \partial_x) + ... \quad ,$$

where $L_l$, $l \geq 1$ are some modified operators of order $2l + 2$. The series does not converge but we consider truncated series:

$$L^{(N)}(\delta, x, \partial_x) = L(x, \partial_x) + \delta L_1(x, \partial_x) + \delta^2 L_2(x, \partial_x) + ... + \delta^N L_N(x, \partial_x) \quad .$$

Note that this operator is no longer of order 2 and we can not define easily a solution to the modified equation

$$\partial_t v^N(t, x) = L^{(N)}(\delta, x, \partial_x)v^N(t, x).$$

However, in our case, we can build an approximated solution $v^{(N)}$ such that

$$\| \mathbb{E}(\phi(X_p)) - v^{(N)}(p\delta, .) \|_{\mathcal{C}} \leq c_2(\phi, N)\delta^N, \quad p = 0, ..., [T/\delta], \quad T > 0,$$

where $\mathcal{C}$ is an appropriate space. As the constant $c_2$ does not depend of $T$, we have an approximation result valid on very long times. We also show that there exists a modified invariant measure for $L^{(N)}(\delta, x, \partial_x)$.

The two main tools are the exponential convergence to equilibrium and the ellipticity of the Poisson equation, i.e. the equation $L(x, \partial_x)u = h$. The second tool is also used in [16]. In a forthcoming article, we will treat the Langevin equation as in [25].

In section 2, we introduce the SDE and the assumptions that we need. We also introduce the numerical schemes that we use. Then, we give some results on the Kolmogorov operator and on the solution of the Kolmogorov equation. In section 3, we give an asymptotic expansion of the weak error. In section 4, we study the modified operator $\mathcal{L}$ and its approximation. In section 5, we analyze the long time behavior of $v^{(N)}$.

## 2 Preliminaries

### 2.1 Presentation of the SDE

In all the article, we write the dot product of two vectors $x = (x_1, ..., x_d) \in \mathbb{R}^d$ and $y = (y_1, ..., y_d) \in \mathbb{R}^d$ as

$$\langle x, y \rangle = \sum_{i=1}^{d} x_i y_i.$$

We identify the gradient of a function with its differential $D$. For $k \in \mathbb{N}$, $h = (h_1, ..., h_k) \in \mathbb{R}^{kd}$, $x \in \mathbb{R}^d$ and $\phi$ a function, we denote the differential of order $k$ of $\phi$ in the point $x$ and in the direction $(h_1, ..., h_k)$ by $D^k \phi(x) \cdot (h_1, ..., h_k)$. For a multi-index $\mathbf{k} = (k_1, ..., k_d) \in \mathbb{N}^d$, we set $|\mathbf{k}| = k_1 + ... + k_d$ and for a function $\phi \in C^\infty(\mathbb{R}^d)$

$$\partial_{\mathbf{k}} \phi(x) = \frac{\partial^{|\mathbf{k}|} \phi(x)}{\partial^{k_1} x_1 ... \partial^{k_d} x_d}, \quad x \in \mathbb{R}^d.$$

We also use the following notation

$$\begin{aligned}
\mathcal{C}^\infty_{pol}(\mathbb{R}^d) =& \{f \in C^\infty(\mathbb{R}^d) \text{ such that } f \text{ and all its derivatives have polynomial growth}\} \\
=& \{f \in C^\infty(\mathbb{R}^d) \text{ such that for all } \mathbf{k} = (k_1, ..., k_d) \in \mathbb{N}^d, \exists C_{\mathbf{k}}, n_{\mathbf{k}} \text{ such that for all } x \in \mathbb{R}^d, \\
& |\partial_{\mathbf{k}} f(x)| \leq C_{\mathbf{k}} (1 + |x|^{2n_{\mathbf{k}}})\}.
\end{aligned}$$

Let $(\Omega, \mathcal{F}, \mathcal{F}_t, \mathbf{P})$, $t \geq 0$, be a filtered probability space and $W(t) = (W_1(t), ..., W_d(t))$ be a $d$-dimensional $\{\mathcal{F}_t\}_{t \geq 0}$-adapted standard Wiener process. We want to give a similar result to the result describe in [4] for a process $X(t)$ on $\mathbb{R}^d$ which verifies the stochastic differential equation

$$dX(t) = -DV(X(t))dt + dW(t) \quad t > 0, \tag{2.1}$$

where $V : \mathbb{R}^d \to \mathbb{R}$, $C^\infty$, verifies the following conditions:

**B-1:** For all $k \in \mathbb{N}$,

$$\int_{\mathbb{R}^d} |x|^{2k} e^{-2V(x)} dx < \infty.$$

**B-2:** The function $V$ is semi-convex : There exist a bounded function $V_1 \in C^\infty(\mathbb{R}^d)$ with bounded derivatives and a convex function $V_2 \in C^\infty(\mathbb{R}^d)$ such that $V = V_1 + V_2$.

**B-3:** There exist strictly positive real numbers $\beta$ and $\kappa$ such that

$$\langle x, DV(x) \rangle \geq \beta |x|^2 - \kappa \quad \text{for all } x \in \mathbb{R}^d.$$

**B-4:** $V \in \mathcal{C}^\infty_{pol}(\mathbb{R}^d)$.

**Remark 2.1** *Under these assumptions, we have that $X(t)$ is well-defined for all $t > 0$ (see [2, 21]).*

A consequence of the semi-convexity of $V$ (assumption B-2) is that there exists a positive constant $\alpha$ such that for all $x \in \mathbb{R}^d$ and $h \in \mathbb{R}^d$

$$D^2 V(x) \cdot (h, h) \geq -\alpha |h|^2. \tag{2.2}$$

In the following, $\alpha$ will be called the constant of semi-convexity of $V$.

We denote by $L$ the Kolmogorov generator associated with the stochastic equation (2.1) : for $\phi \in C^\infty$ and $x \in \mathbb{R}^d$

$$L(x)\phi(x) = \frac{1}{2} \sum_{i=1}^d \partial_{ii} \phi(x) - \sum_{i=1}^d \partial_i V(x) \partial_i \phi(x),$$

where we use the notation $\partial_i = \frac{\partial}{\partial x_i}$ and $\partial_{ij} = \frac{\partial^2}{\partial x_i \partial x_j}$ for $i, j \in \{1, ..., d\}$.

Under assumptions B, we have the following result:

**Proposition 2.2** *Let $x_0 \in \mathbb{R}^d$ and $(X(t))_{t \geq 0}$ satisfying (2.1) and $X(0) = x_0$. Under assumption B-3, we have, for each $p \geq 1$ and $0 < \gamma < 2\beta$, there exists a positive constant $C_p$ such that*

$$\forall t > 0, \quad \mathbb{E}|X(t)|^{2p} \leq C_p \big(|x_0|^{2p} \exp(-\gamma t) + 1\big) \tag{2.3}$$

**Proof.** Let $N \in \mathbb{N}^*$ be fixed. Let $x_0 \in \mathbb{R}^d$ be fixed such that $X(0) = x_0$. We consider $\tau_N = \inf\{t, \text{ such that } |X(t)| \geq N\}$.
Using dissipativity assumption B-3, we get for $x \in \mathbb{R}^d$

$$L|x|^2 = -2\langle x, DV(x)\rangle + d \leq -2\beta|x|^2 + 2\kappa + d. \tag{2.4}$$

We prove (2.3) by recursion. In same time, we will show by recursion the following result: For $p \in \mathbb{N}^*$ and $0 < \gamma < 2\beta$, there exists a positive constant $C_p$ such that for all $t \geq 0$

$$\mathbb{E}\Big(\int_0^{t \wedge \tau_N} |X(s)|^{2p} \exp(\gamma s) ds\Big) \leq C_p\Big(|x_0|^{2p} + 1 + \mathbb{E}\big(\exp(\gamma(t \wedge \tau_N))\big)\Big). \tag{2.5}$$

Let $0 < \gamma_1 < 2\beta$. We apply the Itô's lemma to $|X(t)|^2 \exp(\gamma_1 t)$. We obtain for all $t \geq 0$

$$|X(t \wedge \tau_N)|^2 \exp(\gamma_1(t \wedge \tau_N)) = |X(0)|^2 + \gamma_1 \int_0^{t \wedge \tau_N} |X(s)|^2 \exp(\gamma_1 s) ds$$
$$+ \int_0^{t \wedge \tau_N} L\big(|X(s)|^2\big) \exp(\gamma_1 s) ds + \int_0^{t \wedge \tau_N} 2X(s) \exp(\gamma_1 s) dW(s)$$

because the Brownian motions $W^1, ..., W^d$ are independent.
The stochastic integral is a square integrable martingale because $X(.)$ is bounded on $[0, t \wedge \tau_N]$. Thus, its average vanishes. Using (2.4), we get for all $t \geq 0$

$$\mathbb{E}\Big(|X(t \wedge \tau_N)|^2 \exp(\gamma_1(t \wedge \tau_N))\Big) \leq |x_0|^2 + (\gamma_1 - 2\beta)\mathbb{E}\Big(\int_0^{t \wedge \tau_N} |X(s)|^2 \exp(\gamma_1 s) ds\Big)$$
$$+ \mathbb{E}\Big(\int_0^{t \wedge \tau_N} (2\kappa + d) \exp(\gamma_1 s) ds\Big). \tag{2.6}$$

Using $\gamma_1 < 2\beta$, we get for all $t \geq 0$

$$\mathbb{E}\big(|X(t \wedge \tau_N)|^2 \exp(\gamma_1(t \wedge \tau_N))\big) \leq |x_0|^2 + \frac{2\kappa + d}{\gamma_1}\mathbb{E}\big(\exp(\gamma_1(t \wedge \tau_N))\big).$$

Using Fatou's lemma in the left hand side and Monotone convergence theorem in the right hand side, we obtain for $t \geq 0$

$$\mathbb{E}\big(|X(t)|^2\big)\exp(\gamma_1 t) \leq |x_0|^2 + \frac{2\kappa + d}{\gamma_1}\exp(\gamma_1 t).$$

Thus, the result (2.3) is proved for $p = 1$. Using (2.6) and $\gamma_1 < 2\beta$, we get for $t \geq 0$

$$(2\beta - \gamma_1)\mathbb{E}\Big(\int_0^{t \wedge \tau_N} \big(|X(s)|^2\big)\exp(\gamma_1 s)ds\Big) \leq |x_0|^2 + \frac{2\kappa + d}{\gamma_1}\mathbb{E}\big(\exp(\gamma_1(t \wedge \tau_N))\big).$$

Thus, the result (2.5) is proved for $p = 1$.

Let us assume that (2.3) and (2.5) are true for $p - 1$. We want to obtain them for $p$. We use the same ideas as for $p = 1$.

Let $0 < \gamma_1 < 2\beta$. We apply the Itô's lemma to $\exp(\gamma_1 t)|X(t)|^{2p}$. The average of the stochastic integral vanishes. Using (2.4) and for all $x \in \mathbb{R}^d$ and $h \in \mathbb{R}^d$,

$$D|x|^{2p} \cdot h = 2p\langle x, h\rangle |x|^{2(p-1)}$$

and

$$D^2|x|^{2p} \cdot (h, h) = 2|x|^{2(p-1)}|h|^2 + 4p(p-1)|\langle x, h\rangle|^2|x|^{2(p-2)},$$

we get for $t \geq 0$

$$\begin{aligned}\mathbb{E}\Big(\exp(\gamma_1(t \wedge \tau_N))|X(t \wedge \tau_N)|^{2p}\Big) =& \mathbb{E}\big(|X(0)|^{2p}\big) + \gamma_1\mathbb{E}\Big(\int_0^{t \wedge \tau_N}\exp(\gamma_1 s)|X(s)|^{2p}ds\Big) \\ &+ p\mathbb{E}\Big(\int_0^{t \wedge \tau_N}\exp(\gamma_1 s)L\big(|X(s)|^2\big)|X(s)|^{2(p-1)}ds\Big) \\ &+ 2p(p-1)\mathbb{E}\Big(\int_0^{t \wedge \tau_N}\exp(\gamma_1 s)|X(s)|^{2(p-1)}ds\Big) \\ \leq& |x_0|^{2p} + (\gamma_1 - 2\beta p)\mathbb{E}\Big(\int_0^{t \wedge \tau_N}\exp(\gamma_1 s)|X(s)|^{2p}ds\Big) \\ &+ p(2(p-1) + 2\kappa + d)\mathbb{E}\Big(\int_0^{t \wedge \tau_N}\exp(\gamma_1 s)|X(s)|^{2(p-1)}ds\Big).\end{aligned} \quad (2.7)$$

We use the induction hypothesis. Thus, we get for $t \geq 0$

$$\mathbb{E}\Big(\exp(\gamma_1(t \wedge \tau_N))|X(t \wedge \tau_N)|^{2p}\Big)$$
$$\leq |x_0|^{2p} + p(2(p-1) + 2\kappa + d)C_{p-1}\Big(|x_0|^{2(p-1)} + 1 + \mathbb{E}\big(\exp(\gamma_1(t \wedge \tau_N))\big)\Big)$$

Using Fatou's lemma in the left hand side and Monotone convergence theorem in the right hand side, we obtain for $t \geq 0$

$$\mathbb{E}\big(|X(t)|^{2p}\big)\exp(\gamma_1 t) \leq |x_0|^{2p} + p(2(p-1) + 2\kappa + d)C_{p-1}\Big(|x_0|^{2(p-1)} + 1 + \exp(\gamma_1 t)\Big)$$
$$\leq C_p\Big(|x_0|^{2p} + 1 + \exp(\gamma_1 t)\Big).$$

Thus, we have shown (2.3) for $p$.

Using (2.7) and the induction hypothesis, we get for $t \geq 0$

$$\mathbb{E}\Big( \int_0^{t \wedge \tau_N} |X(s)|^{2p} \exp(\gamma_1 s) ds \Big) \leq C\Big( |x_0|^{2p} + 1 + \mathbb{E}\Big( \exp(\gamma_1(t \wedge \tau_N)) \Big) \Big).$$

We have shown (2.5) for $p$. ♦

## 2.2 Numerical schemes

For a small time step $\delta$ and $x \in \mathbb{R}^d$, the classical explicit Euler method applied to (1.2) is defined for $i = 1, ..., d$, by $X_0 = x$ and the formula

$$X_{n+1}^i = X_n^i + \delta f^i(X_n) + \sum_{l=1}^m g_l^i(X_n)(W^l((n+1)\delta) - W^l(n\delta)), \quad n \geq 0. \tag{2.8}$$

The ordinary Euler scheme (2.8) may be unstable when the coefficients of the differential equation (2.1) are unbounded (see [11]). We are led to avoid explicit schemes. In fact, we study two different implicit schemes. For a small time step $\delta > 0$ and $x \in \mathbb{R}^d$, we consider an implicit split-step scheme defined by $X_0 = x$ and for $n \in \mathbb{N}$

$$\begin{cases} X_n^* & = X_n - \delta DV(X_n^*), \\ X_{n+1} & = X_n^* + W((n+1)\delta) - W(n\delta) = X_n^* + \sqrt{\delta}\eta_n, \end{cases} \tag{2.9}$$

where $\eta_n = (\eta_{n,1}, ..., \eta_{n,d})$ is a $\mathbb{R}^d$-valued random variable and $\{\eta_{n,i} : n \in \mathbb{N}, i \in \{1, ..., d\}\}$ is a collection of i.i.d. real-valued random variables satisfying $\eta_{1,1} \sim \mathcal{N}(0,1)$. We also consider the implicit Euler scheme defined by $X_0 = x$ and for $n \in \mathbb{N}$

$$X_{n+1} = X_n - DV(X_{n+1})\delta + \sqrt{\delta}\eta_n, \tag{2.10}$$

where $\eta_n = (\eta_{n,1}, ..., \eta_{n,d})$ is as above.

Using the following Lemma, we get that these two schemes are well-defined for $\delta < \delta_0 = \frac{1}{\alpha}$, where $\alpha$ is the constant of semi-convexity of $V$.

**Lemma 2.3** *Let $x \in \mathbb{R}^d$ and $\delta < \delta_0 = \frac{1}{\alpha}$. There exists a unique $y \in \mathbb{R}^d$ such that $y = x - \delta DV(y)$.*

**Proof.** Let $x \in \mathbb{R}^d$ and $\delta < \delta_0 = \frac{1}{\alpha}$. Let $P : \mathbb{R}^d \to \mathbb{R}^d$ defined for $z \in \mathbb{R}^d$ by $P(z) = z - x + \delta DV(z)$. We have that $P \in C^\infty$. Using dissipativity assumption B-3, we get

$$\langle P(z), z \rangle = |z|^2 - \langle x, z \rangle + \delta \langle DV(z), z \rangle \geq (1 + \beta\delta)|z|^2 - \kappa\delta - \langle x, z \rangle.$$

Thus, we have for $|z|^2$ large enough that $\langle P(z), z \rangle > 0$. Then, using a Corollary of Brower fixed-point Theorem (see for instance [15]), we have that there exists $y \in \mathbb{R}^d$ such that $P(y) = 0$. Therefore, we have shown the existence of $y \in \mathbb{R}^d$ such that $y = x - \delta DV(y)$.
Let us show the uniqueness.
Let $y \in \mathbb{R}^d$ and $z \in \mathbb{R}^d$ such that $P(z) = P(y) = 0$ and $y \neq z$, then

$$y - z = -\delta(DV(y) - DV(z)).$$

Using assumption of semi-convexity (B-2), we get

$$|y - z|^2 = -\delta \int_0^1 D^2V(y + t(z-y)) \cdot (y - z, y - z) dt \leq \delta\alpha |y - z|^2 < |y - z|^2.$$

Then, $y = z$. ♦

**Remark 2.4** *The condition $\delta < \frac{1}{\alpha}$ is useful only for the uniqueness. We have the existence for all $\delta > 0$. Moreover, in the case where $V$ is convex, the inequality on the second derivative (2.2) is true for all $\alpha \geq 0$. Hence, we can show the uniqueness for all $\delta > 0$.*

**Proposition 2.5** *Let $\delta < \delta_0 := \frac{1}{\alpha}$, where $\alpha$ is the constant of semi-convexity of $V$. Under assumptions B, the implicit split-step scheme (2.9) and the implicit Euler scheme (2.10) satisfy:*

$$\forall p \in \mathbb{N}, \quad \exists C_p(\delta_0) \text{ such that } \forall n \in \mathbb{N}, \quad \mathbb{E}(|X_n|^{2p}) < C_p(\delta_0)(1 + |x|^{2p}). \tag{2.11}$$

**Proof.** We begin to show that the implicit split-step scheme has moments of all order. The proof in the case of the implicit Euler scheme (2.10) is similar.

Let $p \in \mathbb{N}^*$ and $n \in \mathbb{N}^*$ be fixed. Using the convexity of the function $x \mapsto |x|^{2p}$, we have for all $x \in \mathbb{R}^d$ and $y \in \mathbb{R}^d$

$$|x + y|^{2p} \geq |x|^{2p} + 2p|x|^{2(p-1)}(x, y)$$

then

$$|X_n|^{2p} = |X_n^* + \delta DV(X_n^*)|^{2p} \geq |X_n^*|^{2p} + 2p\delta\langle X_n^*, DV(X_n^*)\rangle |X_n^*|^{2(p-1)}.$$

Using dissipativity assumption B-3, we get

$$|X_n|^{2p} \geq (1 + 2p\delta\beta)|X_n^*|^{2p} - 2p\kappa\delta|X_n^*|^{2(p-1)}.$$

Let $\ell \in \mathbb{N}^*$ and $\varepsilon > 0$, then there exists $C_\varepsilon > 0$ such that for any $x \in \mathbb{R}^d$

$$|x|^{2(\ell-1)} \leq \varepsilon|x|^{2\ell} + C_\varepsilon. \tag{2.12}$$

Using (2.12) for $\ell = p$ and $\varepsilon = \frac{\beta}{2\kappa}$, we get

$$|X_n|^{2p} \geq (1 + p\beta\delta)|X_n^*|^{2p} - \delta C. \tag{2.13}$$

Moreover, we have

$$|X_{n+1}|^{2p} = \left(|X_n^*|^2 + 2\sqrt{\delta}\langle X_n^*, \eta_n\rangle + \delta|\eta_n|^2\right)^p$$

$$= \sum_{i+j+k=p} \frac{2^j p!}{i!j!k!}|X_n^*|^{2i}\langle X_n^*, \eta_n\rangle^j|\eta_n|^{2k}\delta^{k+j/2}$$

$$= |X_n^*|^{2p} + \sum_{\substack{i+2j+k=p \\ i \neq p}} \frac{2^{2j}p!}{i!(2j)!k!}|X_n^*|^{2i}\langle X_n^*, \eta_n\rangle^{2j}|\eta_n|^{2k}\delta^{k+j}$$

$$+ \sum_{\substack{i+2j+1+k=p \\ i \neq p}} \frac{2^{2j+1}p!}{i!(2j+1)!k!}|X_n^*|^{2i}\langle X_n^*, \eta_n\rangle^{2j}|\eta_n|^{2k}\langle X_n^*, \eta_n\rangle\delta^{k+j+1/2}$$

$$:= |X_n^*|^{2p} + A + B.$$

Each term in $B$ is a product of an odd number of $\eta_{n,i}$. Then using that for $i \in \{1, ..., d\}$, $\eta_{n,i}$ is independent with $X_n^*$ and $\eta_{n,j}$ for $i \neq j$ and that the moment of odd order of $\eta_{n,i}$ vanish, we have $\mathbb{E}(B) = 0$.

Let $i, j, k \in \mathbb{N}^*$, if $i + 2j + k = p$ and $i \neq p$ then $i + j < p$ and $j \neq 0$ or $k \neq 0$. Let $\varepsilon > 0$, using properties of $\eta_n$ and (2.12), we also have for $i + 2j + k = p$ and $i \neq p$

$$\mathbb{E}\left(|X_n^*|^{2i}\langle X_n^*, \eta_n\rangle^{2j}|\eta_n|^{2k}\right) \leq C\mathbb{E}\left(|X_n^*|^{2(j+i)}\right)$$

$$\leq \varepsilon\mathbb{E}\left(|X_n^*|^{2p}\right) + C_\varepsilon$$

and
$$\mathbb{E}(A) \leq \delta C_{p,d}(\delta_0)(\varepsilon \mathbb{E}|X_n^*|^{2p} + C_\varepsilon).$$

Then, using (2.13), we get

$$\mathbb{E}|X_{n+1}|^{2p} \leq \mathbb{E}|X_n^*|^{2p} + \delta C_{p,d}(\delta_0)\varepsilon \mathbb{E}|X_n^*|^{2p} + \delta C_{p,d}(\delta_0)C_\varepsilon$$
$$\leq \frac{1+\delta C_{p,d}(\delta_0)\varepsilon}{1+p\delta\beta}\mathbb{E}|X_n|^{2p} + \delta C.$$

Choosing $\varepsilon < \frac{p\beta}{C_{p,d}(\delta_0)}$, we get, by induction on $n$,

$$\mathbb{E}|X_n|^{2p} \leq \left(\frac{1+\delta C_{p,d}(\delta_0)\varepsilon}{1+2p\delta\beta}\right)^n |x|^{2p} + \delta \sum_{i=0}^{n-1}\left(\frac{1+\delta C_{p,d}(\delta_0)\varepsilon}{1+2p\delta\beta}\right)^i C$$
$$\leq |x|^{2p} + C(1+|x|^{2(p-1)})$$
$$\leq C_p(\delta_0)(1+|x|^{2p}).$$

♦

## 2.3 Main result

For $x \in \mathbb{R}^d$, we denote by $(X_x(t))_{t\geq 0}$ a process which verifies (2.1) and has for initial data $X_x(0) = x$. From now on, $(P_t)_{t\geq 0}$ is the transition semigroup associated with the Markov process $(X_x(t))_{t\geq 0}$.

We recall that we denote by $L$ the Kolmogorov generator associated with the stochastic equation (2.1) defined for $\phi \in C^\infty(\mathbb{R}^d)$ by

$$L\phi = \frac{1}{2}\sum_{i=1}^d \partial_{ii}\phi - \sum_{i=1}^d \partial_i V \partial_i \phi = \frac{1}{2}e^{2V}\sum_{i=1}^d \partial_i(e^{-2V}\partial_i\phi). \quad (2.14)$$

Moreover, its formal adjoint in $\mathbb{R}^d$ is given for $\phi \in C^\infty(\mathbb{R}^d)$ by

$$L^\top \phi = \sum_{i=1}^d \partial_{ii} V \phi + \frac{1}{2}\sum_{i=1}^d \partial_{ii}\phi + \sum_{i=1}^d \partial_i V \partial_i \phi.$$

We consider $\rho = \frac{1}{Z}e^{-2V}$ where $Z = \int_{\mathbb{R}^d} e^{-2V(x)}dx$. It is classical to prove that the measure $\rho(x)dx =: d\rho$ is invariant by $P_t$ and that $L^T \rho = 0$.
We also define the formal adjoint $L^*$ of $L$ in $L^2(\rho)$. We have the useful equality:

$$L^* = L$$

We will use the following notation for $k \in \mathbb{N}$ and $l \in \mathbb{N}$

$$\mathcal{C}_k^l(\mathbb{R}^d) := \{\phi \in C^l(\mathbb{R}^d) \text{ such that } \exists C_l \text{ such that for all } x \in \mathbb{R}^d \text{ and } \mathbf{j} \in \mathbb{N}^d, |\mathbf{j}| \leq l,$$
$$|\partial_\mathbf{j}\phi(x)| \leq C_l(1+|x|^k)\}.$$

Let $\psi \in \mathcal{C}_k^l(\mathbb{R}^d)$, we define the following norm

$$\|\psi\|_{l,k} := \sup_{\mathbf{j}=(j_1,\ldots j_d), |\mathbf{j}| \leq k}\left(|\partial_\mathbf{j}\psi(x)|(1+|x|^k)^{-1}\right),$$

and the semi-norm
$$|\psi|_{l,k} := \sup_{\mathbf{j}=(j_1,\ldots j_d), 1\leq |\mathbf{j}|\leq l} \left(|\partial_{\mathbf{j}}\psi(x)|(1+|x|^k)^{-1}\right).$$

We consider functions $\psi$ and $\phi$ which are $C^\infty$ and with compact support. Using (2.14), we have by integration by part:
$$\int \phi L\psi d\rho = -\frac{1}{2}\int e^{-2V}(D\phi, D\psi)dx = \int \psi L\phi d\rho. \tag{2.15}$$

This result can be easily extended to $\psi$ and $\phi$ in $\mathcal{C}_{pol}^\infty(\mathbb{R}^d)$.

We consider a function $\phi \in \mathcal{C}_{pol}^\infty(\mathbb{R}^d)$ and we set, for all $x \in \mathbb{R}^d$ and $t \geq 0$,
$$u(t,x) = \mathbb{E}(\phi(X_x(t))) = P_t\phi(x), \tag{2.16}$$

where $(X_x(t))_{t\geq 0}$ is a process which verifies (2.1) and has for initial data $X_x(0) = x$. This is well defined thanks to Proposition 2.2. We classically have that $u$ is a $C^\infty$ function of $(t,x) \in \mathbb{R}^+ \times \mathbb{R}^d$. Moreover, we can show that $u \in \mathcal{C}_{pol}^\infty(\mathbb{R}^d)$ (see Appendix).

It is well known that $u$ is the unique solution of the Kolmogorov equation (see [8]):
$$\frac{d}{dt}u(t,x) = L(x)u(t,x),\ x \in \mathbb{R}^d,\ t > 0, \quad u(0,x) = \phi(x),\ x \in \mathbb{R}^d. \tag{2.17}$$

Note that we use the standard identification $u(t) = u(t,.)$.

We have the following properties which are necessary to have our result. The first is the existence of a solution to the Poisson equation associated with the operator $L$. Under assumptions B, the following Lemma is true (see e.g. [19]).

**Lemma 2.6** *Let $g \in \mathcal{C}_{pol}^\infty(\mathbb{R}^d)$ such that $\int_{\mathbb{R}^d} g(x)\rho(x)dx = 0$. Then, there exists a unique function $\mu \in \mathcal{C}_{pol}^\infty(\mathbb{R}^d)$ such that*
$$L(x)\mu(x) = g(x) \quad and \quad \int_{\mathbb{R}^d} \mu(x)\rho(x)dx = 0. \tag{2.18}$$

The second property is the exponential convergence to 0 of $u$, the solution of (2.17), and its derivatives in an appropriate space:

**Proposition 2.7** *Let $\phi \in \mathcal{C}_{pol}^\infty(\mathbb{R}^d)$. Let $u$ be the solution of (2.17). Assume that $\int_{\mathbb{R}^d}\phi(x)\rho(x)dx = 0$, then there exists a constant $\lambda$ and for each $k \in \mathbb{N}$, there exist integers $n_k$ and $m_k \geq n_k$ and a positive real number $C_k$ which also depend of $V$, such that $\phi \in \mathcal{C}_{n_k}^k(\mathbb{R}^d)$ and we have the following estimate: for all $t \geq 0$ and $x \in \mathbb{R}^d$*
$$|D^k u(t,x)| \leq C_k \parallel \phi \parallel_{k,n_k} e^{-\lambda t}(1+|x|^{m_k}). \tag{2.19}$$

A proof of this proposition can be found in the appendix.

**Remark 2.8** *Lemma 2.6 is a corollary of Proposition 2.7. Indeed, it is shown that the unique solution of (2.18) is defined by*
$$\mu(x) = \int_0^\infty \mathbb{E}\bigl(g(X_x(t))\bigr)dt.$$

Our main result can be stated as follows:

**Theorem 2.9** *Let $X_p$ be the discrete process defined by the implicit Euler scheme (2.10) or the implicit split-step scheme (2.9). Let $N$ and $n_N$ be fixed. Let $\delta_0 = \frac{1}{\alpha}$ where $\alpha$ is the constant of semi-convexity of $V$. Then, for all $\delta < \delta_0$, there exists a modified function*

$$\mu^N = 1 + \sum_{n=1}^{N} \delta^n \mu_n$$

*such that $\mu^N \in \mathcal{C}_{pol}^\infty(\mathbb{R}^d)$ and*

$$\int_{\mathbb{R}^d} \mu^N(x)\rho(x)dx = 1.$$

*For all function $\phi \in \mathcal{C}_{pol}^\infty(\mathbb{R}^d) \cap \mathcal{C}_{n_N}^{6N+2}(\mathbb{R}^d)$, there exist constants $C_N$ and $k_N$ and a positive polynomial function $P_N$ satisfying the following : For all $p \in \mathbb{N}$,*

$$\| \mathbb{E}\phi(X_p) - \int_{\mathbb{R}^d} \phi(x)\mu^N(x)\rho(x)dx \|_{0,k_N} \leq \left( e^{-\lambda t_p} P_N(t_p) + C_N \delta^N \right) \| \phi - \langle \phi \rangle \|_{6N+2,n_N},$$

*where $t_p = p\delta$.*

This result can be viewed as a discrete version of the Proposition 2.7 in the case $k = 0$. We have for $X_p$ the discrete process defined by (2.10) or (2.9) that $\mathbb{E}\phi(X_p)$ which is an approximation of $u$, has the same property as $u$ : $\mathbb{E}\phi(X_p)$ converge exponentially fast to a constant in $\mathcal{C}_{k_N}^0(\mathbb{R}^d)$ up to an error $\delta^N C_N$. At $p$ and $\delta < \delta_0 := \frac{1}{\alpha}$, where $\alpha$ is the constant of semi-convexity of $V$, fixed, we can optimize this error with a good choice of $N$.

Our result can be compared with [4, 16, 23, 24]. As in [4, 16], the only assumption made on $\delta$ is that $\delta < \delta_0$. In the case of $V$ is convex, we do not need to choose $\delta$ smaller than $\delta_0$, we only need that $\delta < \delta_1$ where $\delta_1$ is any fixed number. We also recover an expansion of the invariant measure as in [26]. Our result is similar to the result in the case of SDE on the torus describes in [4] but we need another proofs because all our functions are not bounded.

The constant $C_N$ appearing in the estimate depends of $N$, the semi-convexity constant $\alpha$, the polynomial growth of $V$ and all its derivatives.

## 3 Asymptotic expansion of the weak error

We have the formal expansion for small $t$ and $x \in \mathbb{R}^d$:

$$u(t,x) = \phi(x) + tL(x)\phi(x) + \frac{t^2}{2}L^2(x)\phi(x) + ... + \frac{t^n}{n!}L^n(x)\phi(x) + ... .$$

This is just obtained by Taylor expansion in time.

Since the solution $u(t)$ of the Kolmogorov equation is in $\mathcal{C}_{pol}^\infty(\mathbb{R}^d)$, the above formal expansion can be justified in $L^2(\rho) = \{f : \mathbb{R}^d \to \mathbb{R}, \int |f|^2 d\rho < \infty\}$. Indeed, we have the following easy result whose proof is left to the reader.

**Proposition 3.1** *Let $\delta_1 > 0$ and $\phi \in \mathcal{C}_{pol}^\infty(\mathbb{R}^d)$. Then, for all $N$, there exist a constant $C(N,\phi)$ and an integer $n_1$ which depends of $N$ and of the polynomial growth of $V$, $\phi$ and their derivatives such that for all $\delta < \delta_1$,*

$$|u(\delta,x) - \sum_{n=0}^{N} \frac{\delta^n}{n!} L^n(x)\phi(x)| \leq C(N,\phi)\delta^{N+1}(1 + |x|^{n_1}).$$

We now examine in detail the first time step and its approximation properties in terms of law. By Markov property, it is sufficient to then obtain information at all steps. We want to have an expansion similar to the last Proposition for the process defined by the implicit Euler scheme (2.10) or the implicit split-step scheme (2.9).

**Proposition 3.2** *Let $\phi \in \mathcal{C}_{pol}^\infty(\mathbb{R}^d)$. For any $N \in \mathbb{N}$, there exists an integer $\ell_{2N+2}$ such that $\phi \in \mathcal{C}_{\ell_{2N+2}}^{2N+2}(\mathbb{R}^d)$. For all $n \geq 1$, there exist operators $A_n$ of order $2n$ with coefficients $\mathcal{C}_{pol}^\infty(\mathbb{R}^d)$ which depend of the scheme chosen ((2.10) or (2.9)), such that for all integer $N \geq 1$, there exist a constant $C_N$ and an integer $k$ which depends of $N$, $n_0$ and the polynomial growth of $V$ and its derivatives such that $\forall \delta < \delta_0 := \frac{1}{\alpha}$, where $\alpha$ is the semi-convexity constant of $V$,*

$$|\mathbb{E}\phi(X_1) - \sum_{n=0}^N \delta^n A_n(x)\phi(x)| \leq C_N \delta^{N+1}(1+|x|^k)|\phi|_{2N+2,\ell_{2N+2}}. \tag{3.1}$$

*Moreover $A_0 = I$ and $A_1 = L$.*

**Remark 3.3** *This result is similar of the asymptotic expansion of the weak error describes in [4], but we can not use the same proof. Indeed, we can not use Itô's lemma because the schemes considered here are implicit.*

First, we consider the implicit split-step scheme (2.9). Let $\delta < \frac{1}{\alpha}$ be fixed and $x \in \mathbb{R}^d$ be fixed such that $X_0 = X(0) = x$.
Before proving Proposition 3.2, we need an asymptotic expansion for $X_0^* = x - \delta DV(X_0^*)$. We define the function $\Psi_\delta$ which associate to $z \in \mathbb{R}^d$ the solution $y \in \mathbb{R}^d$ of $y = z - \delta DV(y)$. The function $\Psi_\delta$ is well defined (see Lemma 2.3). Then, by definition of $x$, we have $X_0^* = \Psi_\delta(x)$. Moreover, we have that $(\delta, z) \mapsto \Psi_\delta(z)$ is $C^\infty$ on $]0, \frac{1}{\alpha}[ \times \mathbb{R}^d$. Indeed, let the function $f$ defined on $\Omega_1 = ]0, \frac{1}{\alpha}[ \times \mathbb{R}^{2d}$ by

$$f(\delta, z, y) = z - \delta DV(y) - y,$$

using the assumption of semi-convexity B-2, we have, for all $(\delta, z, y) \in \Omega_1$, that $D_y f(\delta, z, y)$ is invertible. Then, by implicit function Theorem, we have that $(\delta, z) \mapsto \Psi_\delta(z)$ is $C^\infty$ on a neighborhood of each $(\delta, z) \in ]0, \frac{1}{\alpha}[ \times \mathbb{R}^d$.
We have the following lemma:

**Lemma 3.4** *Let $x \in \mathbb{R}^d$ such that $X_0 = x$. We have, for $\delta < \frac{1}{\alpha}$,*

$$\forall N_0 \in \mathbb{N}, \quad X_0^* = y = \Psi_\delta(x) = x + \sum_{k=1}^{N_0} \delta^k d_k(x) + \delta^{N_0+1} R_{N_0+1}(x,\delta), \tag{3.2}$$

*where $\forall k \geq 1$, $d_k \in \mathcal{C}_{pol}^\infty(\mathbb{R}^d)$ is defined for all $z \in \mathbb{R}^d$ by*

$$d_1(z) = -DV(z), \quad \forall k \geq 2, \, d_k(z) = -\sum_{i=1}^{k-1} \frac{1}{i!} \sum_{\substack{k_1+\ldots+k_i=k-1, \\ k_j \geq 1}} D^{i+1}V(z) \cdot (d_{k_1}(z), \ldots, d_{k_i}(z))$$

*and, for any $N \in \mathbb{N}$, $R_{N+1}$ verifies: There exist $C > 0$ and $\ell_N \in \mathbb{N}$ such that for any $z \in \mathbb{R}^d$ and $\delta < \frac{1}{\alpha}$,*

$$|R_{N+1}(z,\delta)| \leq C(1+|z|^{\ell_N}). \tag{3.3}$$

**Proof.** We previously showed that $(\delta, x) \mapsto \Psi_\delta(x)$ is $C^\infty$ on $]0, \frac{1}{\alpha}[\times \mathbb{R}^d$. Then, let $x \in \mathbb{R}^d$ fixed, we have that $d_k(x)$ is the $k^{th}$ term of the Taylor expansion of $\delta \mapsto \Psi_\delta(x)$ and we can write (3.2). We now search an expression for $d_k$.

Let $x \in \mathbb{R}^d$ such that $X_0 = x$ and $y = X_0^*$. Let $N_0$ and $\delta < \frac{1}{\alpha}$ be fixed. We use the temporary notation, for all $1 \leq k \leq N_0$, $g_k = d_k$ and $g_{N_0+1} = R_{N_0+1}$. Using Taylor expansion and

$$y = x + \delta R_1(x, \delta) = x + \sum_{k=1}^{N_0+1} \delta^k g_k(x),$$

we obtain

$$DV(y) = DV(x + \delta R_1(x, \delta))$$
$$= DV(x) + \sum_{n=1}^{N_0-1} \frac{1}{n!} D^{n+1}V(x) \cdot \big(\delta R_1(x, \delta), ..., \delta R_1(x, \delta)\big) + \Theta_{N_0}(x)$$
$$= DV(x) + \sum_{n=1}^{N_0-1} \frac{1}{n!} D^{n+1}V(x) \cdot \bigg(\sum_{k=1}^{N_0+1} \delta^k g_k(x), ..., \sum_{k=1}^{N_0+1} \delta^k g_k(x)\bigg) + \Theta_{N_0}(x)$$
$$= DV(x) + \sum_{n=1}^{N_0-1} \frac{1}{n!} \sum_{m=n}^{n(N_0+1)} \delta^m \sum_{\substack{k_1+...+k_n=m,\\ 1 \leq k_i \leq N_0+1}} D^{n+1}V(x) \cdot \big(g_{k_1}(x), ..., g_{k_n}(x)\big) + \Theta_{N_0}(x)$$
$$= DV(x) + I_1(x) + I_2(x) + \Theta_{N_0}(x),$$

where

$$\Theta_{N_0}(x) = \delta^{N_0} \int_0^1 \frac{(1-t)^{N_0-1}}{(N_0-1)!} D^{N_0+1}V(x + t\delta R_1(x, \delta)) \cdot (R_1(x, \delta), ..., R_1(x, \delta)) dt,$$

$$I_1(x) = \sum_{n=1}^{N_0-1} \frac{1}{n!} \sum_{m=n}^{N_0-1} \delta^m \sum_{\substack{k_1+...+k_n=m,\\ 1 \leq k_i \leq N_0}} D^{n+1}V(x) \cdot (d_{k_1}(x), ..., d_{k_n}(x))$$
$$= \sum_{m=1}^{N_0-1} \delta^m \sum_{n=1}^{m} \frac{1}{n!} \sum_{\substack{k_1+...+k_n=m\\ 1 \leq k_i \leq N_0}} D^{n+1}V(x) \cdot (d_{k_1}(x), ..., d_{k_n}(x)),$$

$$I_2(x) = \sum_{n=1}^{N_0-1} \frac{1}{n!} \sum_{m=N_0}^{n(N_0+1)} \delta^m \sum_{\substack{k_1+...+k_n=m,\\ 1 \leq k_i \leq N_0+1}} D^{n+1}V(x) \cdot (g_{k_1}(x), ..., g_{k_n}(x)).$$

Hence, we get

$$DV(y) = DV(x) + I_1(x) + \delta^{N_0} g(\delta, x),$$

and

$$y = x - \delta DV(y) = x - \delta DV(x) - \delta I_1(x) - \delta^{N_0+1} g(\delta, x), \qquad (3.4)$$

where

$$g(\delta, x) = \int_0^1 \frac{(1-t)^{N_0-1}}{(N_0-1)!} D^{N_0+1}V(x + t\delta R_1(x, \delta)) \cdot (R_1(x, \delta), ..., R_1(x, \delta)) dt$$
$$+ \sum_{n=1}^{N_0-1} \frac{1}{n!} \sum_{m=N_0}^{n(N_0+1)} \delta^{m-N_0} \sum_{\substack{k_1+...+k_n=m,\\ 1 \leq k_i \leq N_0+1}} D^{n+1}V(x) \cdot (g_{k_1}(x), ..., g_{k_n}(x)).$$

Thus, by identifying (3.2) and (3.4) and using the expression of $I_1$ above, we obtain for all $x \in \mathbb{R}^d$,

$$d_1(x) = -DV(x), \quad \forall k \geq 2, \; d_k(x) = -\sum_{i=1}^{k-1} \frac{1}{i!} \sum_{k_1+\ldots+k_i=k-1, k_j \geq 1} D^{i+1}V(x) \cdot (d_{k_1}(x), \ldots, d_{k_i}(x)).$$

Moreover, by induction, we have for all $k \geq 1$, $d_k \in \mathcal{C}^\infty_{pol}(\mathbb{R}^d)$ and $\int d_k d\rho < \infty$. The above identifying does not give an easy expression of $R_{N_0+1}$, then we have not immediately (3.3). To show this result, we will use that, for $N$ fixed, $R_N$ is the remainder of order $N$ of $\delta \mapsto \Psi_\delta(x)$. Therefore, if we show that for $n \in \mathbb{N}$, there exist $C > 0$ and $p_n \in \mathbb{N}$ such that for $\delta < \frac{1}{\alpha}$ and $x \in \mathbb{R}^d$,

$$|\partial_\delta^n \Psi_\delta(x)|^2 \leq C(1 + |x|^{p_n}), \tag{3.5}$$

then we show (3.3) and Lemma 3.4 is shown.

Let $x \in \mathbb{R}^d$ and $\delta < \frac{1}{\alpha}$ be fixed. Let us show (3.5) by induction on $n$.
We have, by definition of $\Psi_\delta(x)$,

$$x = \Psi_\delta(x) + \delta DV\big(\Psi_\delta(x)\big).$$

Then, multiplying by $\Psi_\delta(x)$ and using dissipativity assumption B-3, we get

$$\langle x, \Psi_\delta(x) \rangle = |\Psi_\delta(x)|^2 + \delta \Big\langle DV\big(\Psi_\delta(x)\big), \Psi_\delta(x) \Big\rangle \geq (1 + \beta\delta)|\Psi_\delta(x)|^2 - \kappa\delta.$$

Using

$$2\langle a, b \rangle \leq \varepsilon |a|^2 + \frac{1}{\varepsilon}|b|^2, \text{ for } a \in \mathbb{R}^d, b \in \mathbb{R}^d, \varepsilon > 0, \tag{3.6}$$

we have

$$|\Psi_\delta(x)|^2 (1 + \beta\delta - \frac{\varepsilon}{2}) \leq \frac{1}{2\varepsilon}|x|^2 + \kappa\delta.$$

Taking $\varepsilon = 1$, this shows (3.5) for $n = 0$ and $p_0 = 2$.

Let us assume the result (3.5) is true for all $j < n$ and let us show it for $n$.
We have

$$\partial_\delta^n \Psi_\delta(x) = -\delta \partial_\delta^n \Big(DV\big(\Psi_\delta(x)\big)\Big) - n\partial_\delta^{n-1}\Big(DV\big(\Psi_\delta(x)\big)\Big)$$
$$= -B_1(x, \delta) - nB_2(x, \delta).$$

Using Faà di Bruno's formula, we get

$$B_1(x, \delta) = \delta \sum \frac{n!}{m_1! m_2!(2!)^{m_2} \ldots m_n!(n!)^{m_n}} D^{m_1+\ldots+m_n+1} V\big(\Psi_\delta(x)\big) \cdot \prod_{j=1}^n \big(\partial_\delta^j \Psi_\delta(x)\big)^{m_j},$$

$$B_2(x, \delta) = \sum \frac{(n-1)!}{m_1! m_2!(2!)^{m_2} \ldots m_{n-1}!((n-1)!)^{m_{n-1}}} D^{m_1+\ldots+m_{n-1}+1} V\big(\Psi_\delta(x)\big) \cdot \prod_{j=1}^{n-1} \big(\partial_\delta^j \Psi_\delta(x)\big)^{m_j},$$

where $m_1 + 2m_2 + \ldots + nm_n = n$ in $B_1$ and $m_1 + 2m_2 + \ldots + (n-1)m_{n-1} = n - 1$ in $B_2$.
Multiplying by $\partial_\delta^n \Psi_\delta(x)$, we get

$$|\partial_\delta^n \Psi_\delta(x)|^2 = -\delta D^2 V(\Psi_\delta(x)) \cdot (\partial_\delta^n \Psi_\delta(x), \partial_\delta^n \Psi_\delta(x)) - \langle B_3(x, \delta), \partial_\delta^n \Psi_\delta(x) \rangle - n\langle B_2(x, \delta), \partial_\delta^n \Psi_\delta(x) \rangle,$$

where

$$B_3(x,\delta) = B_1(x,\delta) - \delta D^2 V(\Psi_\delta(x)) \cdot \partial_\delta^n \Psi_\delta(x)$$
$$= \delta \sum \frac{n!}{m_1! m_2! (2!)^{m_2} ... m_n! (n!)^{m_n}} D^{m_1+...+m_n+1} V\left(\Psi_\delta(x)\right) \cdot \prod_{j=1}^{n} \left(\partial_\delta^j \Psi_\delta(x)\right)^{m_j},$$

where $m_1 + 2m_2 + ... + (n-1)m_{n-1} = n$ and $m_n = 0$ in $B_3$.

Then, in $B_2$ and $B_3$, we have that the order of each derivative of $\Psi_\delta$ is less than $n-1$. Hence, using polynomial growth assumption B-4 and induction assumption, we have that $B_2$ and $B_3$ have polynomial growth in $x$. Using (3.6), we have for $\varepsilon > 0$

$$|\partial_\delta^n \Psi_\delta(x)|^2 + \delta D^2 V(\Psi_\delta(x)) \cdot (\partial_\delta^n \Psi_\delta(x), \partial_\delta^n \Psi_\delta(x)) \leq \frac{1}{2\varepsilon} B(x) + \frac{\varepsilon}{2} |\partial_\delta^n \Psi_\delta(x)|^2,$$

where $B$ has polynomial growth. Using semi-convexity assumption B-2 and choosing $\varepsilon$ small enough, we have (3.5). ♦

**Proof of Proposition 3.2 in the case of the implicit split-step scheme** (2.9). Let $x \in \mathbb{R}^d$ such that $X_0 = x$ and $\delta < \frac{1}{\alpha}$. We define $\varepsilon_\delta$ by $\varepsilon_\delta(x) = \Psi_\delta(x) - x$. Let $N$ fixed. We have

$$X_1 = X_0^* + \sqrt{\delta}\eta_0 = x + \varepsilon_\delta(x) + \sqrt{\delta}\eta_0.$$

Using (3.2) and its proof, we have for each $N_0 \in \mathbb{N}^*$,

$$\varepsilon_\delta(x) = \sum_{j=1}^{N_0} \delta^j d_j(x) + \delta^{N_0+1} R_{N_0+1}(x,\delta) = \delta R_1(x,\delta),$$

where $R_{N_0+1}$ and $R_1$ have polynomial growth. Let $\phi \in \mathcal{C}_{pol}^\infty(\mathbb{R}^d)$ such that $\phi \in \mathcal{C}_{n_0}^{2N+2}(\mathbb{R}^d)$. We have

$$\phi(X_1) = \phi(x + \varepsilon_\delta(x) + \sqrt{\delta}\eta_0) = \sum_{k=0}^{2N+1} \frac{1}{k!} \delta^{k/2} D^k \phi(x + \varepsilon_\delta(x)) \cdot (\eta_0,...,\eta_0)$$
$$+ \int_0^1 \frac{(1-t)^{2N+1}}{(2N+1)!} \delta^{N+1} D^{2N+2} \phi(x + \varepsilon_\delta(x) + t\sqrt{\delta}\eta_0) \cdot (\eta_0,...,\eta_0) dt.$$

Let $\lfloor . \rfloor$ denotes the integer part. Using Taylor expansion on $D^k \phi(x + \varepsilon_\delta(x))$ at the order $N_k = N - \lfloor (k+1)/2 \rfloor$ and the computations done in the proof of the previous Lemma, we obtain

$$\phi(X_1) = I_1(x,\eta_0) + I_2(x,\eta_0) + I_3(x,\eta_0) + I_4(x,\eta_0),$$

where

$$I_1(x, \eta_0) = \sum_{k=0}^{2N+1} \frac{1}{k!} \delta^{k/2} D^k \phi(x) \cdot (\eta_0, ..., \eta_0)$$

$$+ \sum_{k=0}^{2N+1} \frac{1}{k!} \delta^{k/2} \sum_{j=1}^{N_k} \frac{1}{j!} \sum_{m=j}^{N_k} \delta^m \sum_{k_1+...+k_j=m, k_s \geq 1} D^{j+k}\phi(x) \cdot (d_{k_1}(x), ..., d_{k_j}(x), \eta_0, ..., \eta_0),$$

$$I_2(x, \eta_0) = \sum_{k=0}^{2N+1} \frac{1}{k!} \delta^{k/2} \sum_{j=1}^{N_k} \frac{1}{j!} \sum_{m=N_k+1}^{jN_k} \delta^m B_{m,j,k}(x, \eta_0),$$

$$I_3(x, \eta_0) = \sum_{k=0}^{2N+1} \frac{1}{k!} \delta^{k/2} \int_0^1 \frac{(1-t)^{N_k}}{N_k!} D^{N_k+1+k} \phi(x + t\varepsilon_\delta(x)) \cdot (\varepsilon_\delta(x), ..., \varepsilon_\delta(x), \eta_0, ...\eta_0) dt,$$

$$I_4(x, \eta_0) = \delta^{N+1} \int_0^1 \frac{(1-t)^{2N+1}}{(2N+1)!} D^{2N+2} \phi(x + \varepsilon_\delta(x) + t\sqrt{\delta}\eta_0) \cdot (\eta_0, ..., \eta_0) dt,$$

and, with the temporary notation, for all $0 \leq k \leq 2N+1$ and $1 \leq i \leq N_k$, $g_{k,i} = d_i$ and $g_{k,N_k+1} = R_{N_k+1}$,

$$B_{m,j,k}(x, \eta_0) = \sum_{k_1+...+k_j=m, k_s \geq 1} D^{j+k}\phi(x) \cdot (g_{k,k_1}(x), ..., g_{k,k_j}(x), \eta_0, ..., \eta_0).$$

We have, for all $p \in \mathbb{N}$ and $i \in \{1, ..., d\}$, $\mathbb{E}(\eta_{0,i}^{2p+1}) = 0$ and $\eta_{0,i}$ is independent with $x$ and $\eta_{0,j}$ for $j \neq i$. Then, the expectation of all the odd term in $k$ in $I_1$, $I_2$ and $I_3$ vanish. Hence, we have

$$\mathbb{E}(I_1(x, \eta_0)) = \sum_{k=0}^{N} \frac{1}{(2k)!} \delta^k \mathbb{E}\Big(D^{2k}\phi(x) \cdot (\eta_0, ..., \eta_0)\Big)$$

$$+ \sum_{k=0}^{N} \frac{1}{(2k)!} \delta^k \sum_{j=1}^{N-k} \frac{1}{j!} \sum_{m=j}^{N-k} \delta^m \sum_{k_1+...+k_j=m, k_s \geq 1} \mathbb{E}\Big(D^{j+2k}\phi(x) \cdot (d_{k_1}(x), ..., d_{k_j}(x), \eta_0, ..., \eta_0)\Big),$$

$$\mathbb{E}(I_2(x, \eta_0)) = \sum_{k=0}^{N} \frac{1}{(2k)!} \delta^k \sum_{j=1}^{N-k} \frac{1}{j!} \sum_{m=N+1-k}^{j(N-k)} \delta^m \mathbb{E}(B_{m,j,2k}(x, \eta_0))$$

$$= \delta^{N+1} \sum_{k=0}^{N} \frac{1}{(2k)!} \sum_{j=1}^{N-k} \frac{1}{j!} \sum_{m=N+1-k}^{j(N-k)} \delta^{m+k-N-1} \mathbb{E}(B_{m,j,2k}(x, \eta_0)),$$

$$\mathbb{E}(I_3(x, \eta_0)) = \sum_{k=0}^{N} \frac{1}{(2k)!} \delta^k \int_0^1 \frac{(1-t)^{N-k}}{(N-k)!} \mathbb{E}\Big(D^{N+1+k}\phi(x + t\varepsilon_\delta(x)) \cdot (\varepsilon_\delta(x), ..., \varepsilon_\delta(x), \eta_0, ...\eta_0)\Big) dt,$$

$$\mathbb{E}(I_4(x, \eta_0)) = \delta^{N+1} \int_0^1 \frac{(1-t)^{2N+1}}{(2N+1)!} \mathbb{E}\Big(D^{2N+2}\phi(x + \varepsilon_\delta(x) + t\sqrt{\delta}\eta_0) \cdot (\eta_0, ..., \eta_0)\Big) dt.$$

Using the polynomial growth of $d_i$ and $R_i$ for all $i$ and $\phi \in \mathcal{C}^{2N+2}_{\ell_{2N+2}}(\mathbb{R}^d)$, we have that there exists an integer $n_1$ which depends of $N$, $\ell_{2N+2}$ and the polynomial growth of $V$ and its derivatives such that

$$|\mathbb{E}(I_2(x, \eta_0))| \leq C_N \delta^{N+1}(1 + |x|^{n_1})|\phi|_{2N,\ell_{2N+2}}.$$

Similarly, we have that there exists an integer $n_2$ which depends of $N$, $\ell_{2N+2}$ and the polynomial growth of $V$ and its derivatives such that

$$|\mathbb{E}(I_4(x, \eta_0))| \leq C_N \delta^{N+1}(1 + |x|^{n_2}) \| D^{2N+2}\phi \|_{0,\ell_{2N+2}}.$$

Using $\varepsilon_\delta(x) = \delta R_1(x,\delta)$ and similar computations, we have that there exists an integer $n_3$ which depends of $N$, $\ell_{2N+2}$ and the polynomial growth of $V$ and its derivatives such that

$$|\mathbb{E}(I_3(x,\eta_0))| \leq C_N \delta^{N+1}(1+|x|^{n_3}) \parallel D^{N+1}\phi \parallel_{N,\ell_{2N+2}}.$$

Moreover, we have

$$\mathbb{E}(I_1(x,\eta_0)) = \sum_{n=0}^{N} \delta^n A_n(x)\phi(x),$$

where $A_0 = I$ and, for all $n \in \mathbb{N}^*$ and $x \in \mathbb{R}^d$,

$$A_n(x)\phi(x) = \sum_{\substack{m+k=n,\\ m\geq 1, k\geq 0}} \Big( \sum_{j=1}^{m} \frac{1}{j!(2k)!} \sum_{k_1+\ldots k_j=m, k_s\geq 1} \mathbb{E}\big(D^{j+2k}\phi(x)\cdot(d_{k_1}(x),\ldots,d_{k_j}(x),\eta_0,\ldots,\eta_0))\big)\Big)$$

$$+ \frac{1}{(2n)!}\mathbb{E}(D^{2n}\phi(x)\cdot(\eta_0,\ldots,\eta_0)).$$

For $n=1$, we have only one possibility: $m=1$ and $k=0$, then $j=1$, $k_1=1$ and

$$A_1\phi = \mathbb{E}(D\phi(x)\cdot d_1(x)) + \frac{1}{2}\mathbb{E}(D^2\phi(x)\cdot(\eta_0,\eta_0)).$$

Using properties of $\eta_0$ and the definition of $d_1$, we get $A_1 = L$.

We have that there exists an integer $k$ which depends of $N$, $\ell_{2N+2}$ and the polynomial growth of $V$ and their derivatives such that

$$|\mathbb{E}\phi(X_1) - \sum_{n=0}^{N} \delta^n A_n(x)\phi(x)| \leq C_N \delta^{N+1}(1+|x|^k)|\phi|_{2N+2,\ell_{2N+2}}.$$

Finally, we have proved (3.1) for the implicit split-step scheme (2.9). ♦

The proof in the case of the implicit Euler scheme (2.10) use the same arguments. Let $x \in \mathbb{R}^d$ such that $X_0 = x$ and $\delta < \frac{1}{\alpha}$. We need an asymptotic expansion for $X_1 = x - DV(X_1)\delta + \sqrt{\delta}\eta_0$. We use the local notation $\theta = \sqrt{\delta}$. We define the function $\psi_\theta$ which associate to $y$ the solution $z$ of $z = x - \theta^2 DV(z) + \theta\eta_0$. This function is well defined (see Lemma 2.3) and we have $X_1 = \psi_\theta(x)$. If we consider the function $f_1$ defined on $]0,\frac{1}{\sqrt{\alpha}}[\times\mathbb{R}^d\times\mathbb{R}^d$ by

$$f_1(\theta,y,z) = -z + y - \theta^2 DV(z) + \theta\eta_0,$$

we can show, as previously, that $(\theta,y) \mapsto \psi_\theta(y)$ is $C^\infty$ on $]0,\frac{1}{\sqrt{\alpha}}[\times\mathbb{R}^d$.

We have the following lemma:

**Lemma 3.5** *Let $x \in \mathbb{R}^d$ such that $X_0 = x$. For $\delta < \frac{1}{\alpha}$, with the local notation $\theta = \sqrt{\delta}$*

$$\forall N_0 \in \mathbb{N}, \quad X_1 = \psi_\theta(x) = x + \sum_{k=1}^{N_0} \delta^{\frac{k}{2}} d_k(x,\eta_0) + \delta^{\frac{N_0+1}{2}} R_{N_0}(x,\delta,\eta_0)$$

*where $\forall k \geq 1$, $d_k$ is defined for all $z \in \mathbb{R}^d$ by*

$$d_1(z,\eta_0) = \eta_0, \quad d_2(z,\eta_0) = -DV(z),$$

$$\forall k \geq 3, \ d_k(z,\eta_0) = \sum_{i=1}^{k-2} \frac{1}{i!} \sum_{k_1+\ldots+k_i=k-2, k_j\geq 1} -D^{i+1}V(z)\cdot(d_{k_1}(z,\eta_0),\ldots,d_{k_i}(z,\eta_0)).$$

*Moreover, we have that* $\mathbb{E}(d_k) \in \mathcal{C}_{pol}^{\infty}(\mathbb{R}^d)$, *for all* $z \in \mathbb{R}^d$ *and* $k \in \mathbb{N}$

$$\mathbb{E}(d_{2k+1}(z, \eta_0)) = 0 \tag{3.7}$$

*and, for any* $N \in \mathbb{N}$, $R_{N+1}$ *verifies: There exist* $C > 0$ *and* $\ell_N \in \mathbb{N}$ *such that for any* $z \in \mathbb{R}^d$ *and* $\delta < \frac{1}{\alpha}$,

$$|\mathbb{E}(R_{N+1}(z, \delta, \eta_0))| \leq C(1 + |z|^{\ell_N}).$$

**Proof.** To prove this Lemma, we use the same ideas as in Lemma 3.4. We first compute $d_k$ for all $k$. By induction, we rewrite $d_k$ only in terms of $d_1$, $d_2$ and the derivatives of $V$. Using the independence of $\eta_0$ with $x$, we can show (3.7).

To prove that $\mathbb{E}(R_{N_0})$ has polynomial growth, we show that for any $n \in \mathbb{N}$, there exist $C_n > 0$ and $k_n \in \mathbb{N}$ such that for $x \in \mathbb{R}^d$, $\delta < \frac{1}{\alpha}$ and the local notation $\theta = \sqrt{\delta}$,

$$\mathbb{E}(|\partial_\theta \psi_\theta(x)|^2) \leq C_n(1 + |x|^{k_n}).$$

♦

The proof of Proposition 3.2 in the case of the implicit Euler scheme (2.10) is similar to the case of the implicit split-step scheme (2.9), but we must use an asymptotic expansion of $D^k \phi$ to a larger order ($2N + 1 - k$ instead of $N - \lfloor (k+1)/2 \rfloor$).

## 4 Modified generator

For now on, all definitions depend of the scheme considered.

### 4.1 Formal series analysis

Let us now consider $\delta < \frac{1}{\alpha}$ as fixed. We want to construct a formal series

$$\mathcal{L} = L + \delta L_1 + \ldots + \delta^n L_n + \ldots \tag{4.1}$$

where the coefficients of the operator $L_n$ are in $\in \mathcal{C}_{pol}^{\infty}(\mathbb{R}^d)$, and such that formally the solution $v$ at time $t = \delta$ of the equation

$$\partial_t v(t, x) = \mathcal{L} v(t, x), t > 0, x \in \mathbb{R}^d \quad v(0, x) = \phi(x), x \in \mathbb{R}^d$$

coincides in the sense of asymptotic expansion with the approximation of the transition semigroup $\mathbb{E}\phi(X_1)$ studied in the previous section. In other words, we want to have the equality in the sense of asymptotic expansion in powers of $\delta$

$$\exp(\delta \mathcal{L})\phi = \phi + \sum_{n \geq 1} \delta^n A_n \phi,$$

where the operators $A_n$ are defined in Proposition 3.2.

Formally, this equation can be written as

$$\exp(\delta \mathcal{L}) - I_d = \delta \tilde{A}(\delta), \tag{4.2}$$

where $\tilde{A}(\delta) = \sum_{n \geq 1} \delta^{n-1} A_n$.

We have
$$\exp(\delta\mathcal{L}) - I_d = \delta\mathcal{L}\bigg(\sum_{n\geq 0} \frac{\delta^n}{(n+1)!}\mathcal{L}^n\bigg).$$

Note that the (formal) inverse of the series is given by
$$\bigg(\sum_{n\geq 0} \frac{\delta^n}{(n+1)!}\mathcal{L}^n\bigg)^{-1} = \sum_{n\geq 0} \frac{B_n}{n!}\delta^n\mathcal{L}^n,$$

where the $B_n$ are the Bernoulli numbers (see [6, 10]). Hence, equations (4.1) and (4.2) are equivalent in the sense of formal series to
$$\mathcal{L} = \sum_{\ell\geq 0} \frac{B_\ell}{\ell!}\delta^\ell \mathcal{L}^\ell \tilde{A}(\delta) = \sum_{n\geq 0} \delta^n \bigg( A_{n+1} + \sum_{\ell=1}^{n} \frac{B_\ell}{\ell!} \sum_{n_1+\ldots+n_{\ell+1}=n-\ell} L_{n_1}\ldots L_{n_\ell} A_{n_{\ell+1}+1} \bigg). \tag{4.3}$$

Identifying the right hand sides of (4.1) and (4.3), we get the following induction formula
$$L_n = A_{n+1} + \sum_{\ell=1}^{n} \frac{B_\ell}{\ell!} \sum_{n_1+\ldots+n_{\ell+1}=n-\ell} L_{n_1}\ldots L_{n_\ell} A_{n_{\ell+1}+1}. \tag{4.4}$$

Each of the terms of the above sum is an operator of order $2n+2$ with coefficients $\mathcal{C}_{pol}^\infty(\mathbb{R}^d)$ and therefore $L_n$ is also an operator of order $2n+2$ with coefficients $\mathcal{C}_{pol}^\infty(\mathbb{R}^d)$.

Notes that (4.2) gives immediately the inverse relation of this formal series equation:
$$A_n = \sum_{\ell=1}^{n} \frac{1}{\ell!} \sum_{n_1+\ldots+n_\ell=n-\ell} L_{n_1}\ldots L_{n_\ell}. \tag{4.5}$$

Moreover, we have clearly
$$L_n 1 = 0.$$

## 4.2 Approximate solution of the modified flow

For a given $N$, we have constructed in the previous section a modified operator
$$L^{(N)} = L + \sum_{n=1}^{N} \delta^n L_n. \tag{4.6}$$

In order to perform weak backward error analysis and estimate recursively the modified invariant law of the numerical process, we should be able to define a solution $v^N$ of the modified flow
$$\partial_t v^N(t,x) = L^{(N)} v^N(t,x), t > 0, x \in \mathbb{R}^d \quad v^N(0,x) = \phi(x), x \in \mathbb{R}^d. \tag{4.7}$$

However, in our situation we do not know whether this equation has a solution.

The goal of the following theorem is to give a proper definition of a modified flow associated to (4.6).

**Theorem 4.1** *Let $\phi \in \mathcal{C}_{pol}^\infty(\mathbb{R}^d)$ such that $\int \phi d\rho = 0$. For any $N \in \mathbb{N}$, there exists an integer $\ell_N$ such that $\phi \in \mathcal{C}_{\ell_N}^N(\mathbb{R}^d)$. For all $n \in \mathbb{N}$, there exist functions $v_n(t,.) \in \mathcal{C}_{pol}^\infty(\mathbb{R}^d)$ defined for all times $t \geq 0$ such that for all $t > 0$, $x \in \mathbb{R}^d$ and $n \in \mathbb{N}$,*
$$\partial_t v_n(t,x) - Lv_n(t,x) = \sum_{\ell=1}^{n} L_\ell v_{n-\ell}(t,x), \tag{4.8}$$

with initial condition $v_0(0, x) = \phi(x)$ for all $x \in \mathbb{R}^d$ and $v_n(0, x) = 0$ for $n \geq 1$ and $x \in \mathbb{R}^d$. For all $N \geq 0$, $x \in \mathbb{R}^d$ and $t \geq 0$, setting

$$v^{(N)}(t, x) = \sum_{k=0}^{N} \delta^k v_k(t, x),$$

then the following holds:

**a.** For $\delta_0 = \frac{1}{\alpha}$, there exist constants $C_N$, $k_N$ and $r_N$ such that for all $t \geq 0$, $x \in \mathbb{R}^d$ and $\delta < \delta_0$,

$$|\mathbb{E}v^{(N)}(t, X_1) - v^{(N)}(t + \delta, x)|$$
$$\leq \delta^{N+1} C_N (1 + |x|^{r_N}) \sup_{\substack{s \in ]0, \delta[, \\ n=0,\ldots,N}} |v_n(t + s, .)|_{2N+2, k_N}.$$

**b.** For $\delta_0 = \frac{1}{\alpha}$, there exist constants $C_N$ and $r_N$ such that for all $\delta < \delta_0$,

$$\| \mathbb{E}\phi(X_1) - v^{(N)}(\delta, .) \|_{0, r_N} \leq \delta^{N+1} C_N \| \phi \|_{6N+2, \ell_{6N+2}}.$$

**Proof.** For $n = 0$, equation (4.8) reduces to $v_0 = u$, the solution of (2.17). By Proposition 2.7, we have that $u$ and all its derivatives have polynomial growth in space and exponential decrease in time. Let $n \in \mathbb{N}^*$ and assume that $v_j$ are constructed for $j = 1, \ldots, n-1$. Let for $x \in \mathbb{R}^d$ and $t \geq 0$

$$F_n(t, x) = \sum_{\ell=1}^{n} L_\ell v_{n-\ell}(t, x), \tag{4.9}$$

the right-hand side in (4.8). Then $v_n$ is uniquely defined and given by the formula

$$v_n(t, .) = \int_0^t P_{t-s} F_n(s, .) ds, \quad t \geq 0. \tag{4.10}$$

Using an induction argument and Proposition 2.7, we know that $v_n$ and all its derivatives have polynomial growth. Moreover, we have for all $i \in \mathbb{N}^*$ and $n \in \mathbb{N}^*$ that there exist integer $k_{n,i}$ and $j_{k,n}$ such that for all $t \geq 0$,

$$\| v_n(t) \|_{i, k_{n,i}} \leq P(t) \| \phi \|_{k+4n, j_{k,n}}, \tag{4.11}$$

where $P$ is a polynomial in $t$ which also depends on $k$, $n$ and $V$. This proves the first part of the Theorem.

To prove **a.**, we consider a fixed time t, and define the functions $w_n(s, x) := v_n(t + s, x)$ for $s \geq 0$, $x \in \mathbb{R}^d$ and $n \in \mathbb{N}$. By definition, these functions satisfy the relations

$$\partial_s w_n(s, x) = \sum_{\ell=0}^{n} L_\ell w_{n-\ell}(s, x), s > 0, x \in \mathbb{R}^d \quad w_n(0, x) = v_n(t, x), x \in \mathbb{R}^d.$$

Let us consider the successive time derivatives of the functions $w_n$. We have, using the definition of $w_n$, for all $s > 0$ and $x \in \mathbb{R}^d$

$$\partial_s^2 w_n(s, x) = \sum_{\ell=0}^{n} L_\ell \partial_s w_{n-\ell}(s, x) = \sum_{k=0}^{n} \sum_{\ell_1 + \ell_2 = k} L_{\ell_1} L_{\ell_2} w_{n-k}(s, x),$$

and we see by induction that for all $m \geq 1$, $x \in \mathbb{R}^d$ and $s > 0$

$$\partial_s^m w_n(s,x) = \sum_{\ell_1+...+\ell_{m+1}=n} L_{\ell_1}...L_{\ell_m} w_{\ell_{m+1}}(s,x).$$

Using the fact that the operators $L_\ell$ are of order $2\ell + 2$ with no terms of order zero and the coefficients of $L_\ell$ have polynomial growth, we see that there exist a constant $C$ depending on $n$ and $m$ and a constant $r_n$, such that for all $s > 0$, $x \in \mathbb{R}^d$ and $m \geq 1$

$$|\partial_s^m w_n(s,x)| \leq C(1+|x|^{r_n}) \sup_{\substack{k=0,...,n \\ 1 \leq j \leq 2(n-k)+2m}} |\partial_j w_k(s,x)|.$$

Now let us consider the Taylor expansion of $w_n(\delta,.)$, for $\delta < \delta_0$. We have for $\delta < \delta_0$, $x \in \mathbb{R}^d$ and $n = 0,...,N$,

$$w_n(\delta, x) = \sum_{m=0}^{N-n} \frac{\delta^m}{m!} \partial_s^m w_n(0,x) + \int_0^\delta \frac{\sigma^{N-n}}{(N-n)!} \partial_s^{N-n+1} w_n(\sigma, x) d\sigma$$

$$= \sum_{m=0}^{N-n} \frac{\delta^m}{m!} \sum_{\ell_1+...+\ell_{m+1}=n} L_{\ell_1}...L_{\ell_m} w_{\ell_{m+1}}(0,x) + R_{N,n}(\delta, x).$$

Using the bounds on the time derivatives of $w_n$, we obtain that there exists a constant $\ell_N$ such that for all $0 < \delta < \delta_0$, $x \in \mathbb{R}^d$ and all $n = 0,...,N$,

$$|R_{N,n}(\delta, x)| \leq C\delta^{N-n+1}(1+|x|^{r_N}) \sup_{\substack{s \in ]0,\delta[, \\ i=0,...,N \\ 1 \leq j \leq 2N+2}} |\partial_j w_i(s,x)|$$

for some constants depending on $N$, $n$. After summation in $n$, and using the expression (4.5) of the operators $A_n$ and the definition of $w_n$, we get for all $x \in \mathbb{R}^d$, $t \geq 0$ and $0 < \delta < \delta_0$

$$v^{(N)}(t+\delta, x) = \sum_{n=0}^{N} \delta^n \sum_{m=0}^{n} A_m v_{n-m}(t,x) + R_N(t,\delta,x)$$

where for all $t \geq 0$ and $0 < \delta < \delta_0$,

$$|R_N(t,\delta,x)| \leq C_N \delta^{N+1}(1+|x|^{r_N}) \sup_{\substack{s \in ]0,\delta[, \\ n=0,...,N \\ 1 \leq j \leq 2N+2}} |\partial_j v_n(t+s,x)|.$$

To conclude, we use (3.1) applied to $\phi = v^{(N)}(t)$ and the fact that $\delta < \delta_0$.
The second estimate **b.** is then a consequence of **a.** with $t = 0$ and (4.11).

## 5 Asymptotic expansion of the invariant measure and long time behavior

We now analyze the long time behavior of the solution of the modified flow associated to (4.7). In the following, for a given operator $B$, we denote by $B^*$ its formal adjoint with respect to the $L^2(\rho)$ product. We start by an asymptotic expansion of a formal invariant measure for the numerical schemes.

**Proposition 5.1** Let $\delta_0 < \frac{1}{\alpha}$. Let $(L_n)_{n \geq 0}$ be the collection of operators defined recursively by (4.4). There exists a collection of functions $(\mu_n)_{n \geq 0}$ such that $\mu_0 = 1$, $\int_\mathbb{R} \mu_n(x)\rho(x)dx = 0$ for $n \geq 1$, and for all $n \geq 1$, $\mu_n \in \mathcal{C}^\infty_{pol}(\mathbb{R}^d)$ and

$$L\mu_n = -\sum_{\ell=1}^n (L_\ell)^* \mu_{n-\ell}. \tag{5.1}$$

Let $N \geq 0$ be fixed and the function $\mu^{(N)}$ be defined for $x \in \mathbb{R}^d$ and $0 < \delta < \delta_0$ by

$$\mu^{(N)}(\delta, x) = 1 + \sum_{n=1}^N \delta^n \mu_n(x).$$

Then for all $0 < \delta < \delta_0$, $\mu^{(N)}(\delta, .) \in \mathcal{C}^\infty_{pol}(\mathbb{R}^d)$ and $\mu^{(N)}$ satisfies for $0 < \delta < \delta_0$

$$\int_{\mathbb{R}^d} \mu^{(N)}(\delta, x)\rho(x)dx = 1.$$

**Remark 5.2** We consider equation (5.1) because $L^* = L$.

**Proof.** Let $n \geq 1$. Assume that $\mu_0 = 1$ and $\mu_j$ are known, for $j = 1, ..., n-1$. Let us consider equation (5.1) given by

$$L\mu_n = -\sum_{\ell=1}^n (L_\ell)^* \mu_{n-\ell} =: G_n$$

Note that $G_n \in \mathcal{C}^\infty_{pol}(\mathbb{R}^d)$ and satisfies

$$\int_{\mathbb{R}^d} G_n(x)\rho(x)dx = -\sum_{\ell=1}^n \int_{\mathbb{R}^d} (L_\ell)^* \mu_{n-\ell}(x)\rho(x)dx = -\sum_{\ell=1}^n \int_{\mathbb{R}^d} \mu_{n-\ell}(x) L_\ell 1 \rho(x)dx = 0.$$

Using the Lemma 2.6, we easily obtain the existence of a function $\mu_n \in \mathcal{C}^\infty_{pol}(\mathbb{R}^d)$ satisfying (5.1) and $\int_{\mathbb{R}^d} \mu_n(x)\rho(x)dx = 0$. This shows the proposition. ♦

**Proposition 5.3** Let $\phi \in \mathcal{C}^\infty_{pol}(\mathbb{R}^d)$. For all $n$ and $k$ there exist a positive polynomial function $P_{k,n}$ and integers $\ell_{n,k}$ and $m_{n,k} \geq \ell_{n,k}$ such that $\phi \in \mathcal{C}^{k+4n}_{\ell_{k,n}}(\mathbb{R}^d)$ and for all $t \geq 0$

$$\| v_n(t) - \int_{\mathbb{R}^d} \phi(y)\mu_n(y)\rho(y)dy \|_{k,m_{n,k}} \leq P_{k,n}(t)e^{-\lambda t} \| \phi - \langle\phi\rangle \|_{k+4n, \ell_{k,n}}, \tag{5.2}$$

where $\langle\phi\rangle = \int_{\mathbb{R}^d} \phi(x)\rho(x)dx$.

**Proof.** Using the fact that $\mu_0 = 1$ and $v_0 = u$, we see that estimate (5.2) is satisfied for $n = 0$ (Proposition 2.7). Let $n \geq 1$ and assume that $v_j$, $j = 0, ..., n-1$ satisfy for $k \in \mathbb{N}^*$, $i \in \mathbb{N}^*$ and $t \geq 0$:

$$\| v_j(t) - \int_{\mathbb{R}^d} \phi(y)\mu_j(y)\rho(y)dy \|_{k,m_{j,k}} \leq P_{k,j}(t)e^{-\lambda t} \| \phi - \langle\phi\rangle \|_{k+4j, \ell_{k,j}},$$

where $\ell_{k,j}$ is such that $\phi \in \mathcal{C}^{k+4j}_{\ell_{k,j}}(\mathbb{R}^d)$. Let us set for $t \geq 0$

$$c_n(t) = \sum_{m=0}^n \int_{\mathbb{R}^d} v_{n-m}(t, x)\mu_m(x)\rho(x)dx.$$

We claim that $c_n$ does not depend on time. Indeed, for all $t \geq 0$,

$$\sum_{m=0}^{n} \partial_t \int_{\mathbb{R}^d} v_{n-m}(t,x) \mu_m(x) \rho(x) dx = \sum_{m=0}^{n} \partial_t \int_{\mathbb{R}^d} v_m(t,x) \mu_{n-m}(x) \rho(x) dx$$

$$= \sum_{m=0}^{n} \sum_{\ell=0}^{m} \int_{\mathbb{R}^d} L_{m-\ell} v_\ell(t,x) \mu_{n-m}(x) \rho(x) dx$$

$$= \sum_{\ell=0}^{n-1} \int_{\mathbb{R}^d} v_\ell(t,x) \sum_{m=1}^{n-\ell} L_m^* \mu_{n-\ell-m}(x) \rho(x) dx$$

$$+ \sum_{\ell=0}^{n-1} \int_{\mathbb{R}^d} v_\ell(t,x) L \mu_{n-\ell}(x) \rho(x) dx + \int_{\mathbb{R}^d} v_n(t,x) L 1 d\rho$$

$$= 0,$$

by definition of the coefficients $\mu_n$ (see (5.1)) and by (2.15). Note that the computation above is justified because $\forall n$, $v_n$ and $\mu_n$ are in $\mathcal{C}_{pol}^\infty(\mathbb{R}^d)$. We deduce for all $t \geq 0$

$$\int_{\mathbb{R}^d} \partial_t v_n(t,x) \rho(x) dx = -\sum_{m=1}^{n} \int_{\mathbb{R}^d} \partial_t v_{n-m}(t,x) \mu_m(x) \rho(x) dx. \tag{5.3}$$

Now, we compute the average of $F_n$. By (4.8), (4.9) and (5.3), we have for all $t \geq 0$

$$\langle F_n(t) \rangle = \int_{\mathbb{R}^d} F_n(t,x) \rho(x) dx = \int_{\mathbb{R}^d} \partial_t v_n(t,x) \rho(x) dx - \int_{\mathbb{R}^d} L v_n(t,x) \rho(x) dx$$

$$= \int_{\mathbb{R}^d} \partial_t v_n(t,x) \rho(x) dx$$

$$= -\sum_{m=1}^{n} \int_{\mathbb{R}^d} \partial_t v_{n-m}(t,x) \mu_m(x) dx.$$

We rewrite (4.10) as follows : for all $t \geq 0$ and $x \in \mathbb{R}^d$

$$v_n(t,x) = \int_0^t \langle F_n(s) \rangle ds + \int_0^t P_{t-s}(F_n(s,x) - \langle F_n(s) \rangle) ds.$$

Using the previous expression obtained for $\langle F_n(s) \rangle$ and recalling the initial data for $v_n$, we deduce that for all $x \in \mathbb{R}^d$ and $t \geq 0$

$$v_n(t,x) = -\sum_{m=1}^{n} \int_{\mathbb{R}^d} v_{n-m}(t,y) \mu_m(y) \rho(y) dy + \int_{\mathbb{R}^d} \phi(y) \mu_n(y) \rho(y) dy + \int_0^t P_{t-s}(F_n(s,x) - \langle F_n(s) \rangle) ds.$$

Then, using $\int_{\mathbb{R}^d} \mu_m(x) \rho(x) dx = 0$, for $m \in \mathbb{N}^*$ (Proposition 5.1), we get for $x \in \mathbb{R}^d$ and $t \geq 0$

$$v_n(t,x) - \int_{\mathbb{R}^d} \phi(y) \mu_n(y) \rho(y) dy = -\sum_{m=1}^{n} \int_{\mathbb{R}^d} \left( v_{n-m}(t,y) - \int_{\mathbb{R}^d} \phi(z) \mu_{n-m}(z) \rho(z) dz \right) \mu_m(y) \rho(y) dy$$

$$+ \int_0^t P_{t-s}(F_n(s,x) - \langle F_n(s) \rangle) ds.$$

Note that, since $L_\ell$, $\ell \in \mathbb{N}$ is a differential operator of order $2\ell+2$ whose the coefficients belong to $\mathcal{C}^\infty_{pol}(\mathbb{R}^d)$ and contain no zero order terms, then we have that there exists an integer $\beta_{n,k}$ such that for $s \geq 0$

$$\| F_n(s) - \langle F_n(s)\rangle \|_{k,\beta_{n,k}} \leq \sum_{\ell=0}^{n-1} c_{k,\ell} |v_\ell(s)|_{2(n-\ell)+2+k, m_{\ell,k+2(n-\ell)+2}}$$

$$\leq \sum_{\ell=0}^{n-1} c_{k,\ell} \| v_\ell(s) - \int_{\mathbb{R}^d} \phi(y)\mu_\ell(y)\rho(y)dy \|_{2(n-\ell)+2+k, m_{\ell,k+2(n-\ell)+2}}.$$

We have used

$$|v_\ell(s)|_{2(n-\ell)+2+k, m_{\ell,k+2(n-\ell)+2}} = |v_\ell(s) - \int_{\mathbb{R}^d} \phi(y)\mu_\ell(y)\rho(y)dy|_{2(n-\ell)+2+k, m_{\ell,k+2(n-\ell)+2}}.$$

Moreover, using the Proposition 2.7, we have for $t \geq 0$, $i \in \mathbb{N}^*$, $x \in \mathbb{R}^d$ and $\mathbf{k} \in \mathbb{N}^d$ that there exists a real number $\alpha_{\mathbf{k},n}$ such that

$$|\partial_\mathbf{k}(v_n(t,x) - \int_{\mathbb{R}^d} \phi(y)\mu_n(y)\rho(y)dy)| \leq$$

$$\sum_{m=1}^n \Big( \int_{\mathbb{R}^d} |v_{n-m}(t,y) - \int_{\mathbb{R}^d} \phi(z)\mu_{n-m}(z)\rho(z)dz|^2 \rho(y)dy \Big)^{1/2} \Big( \int_{\mathbb{R}^d} |\mu_m(y)|^2 \rho(y)dy \Big)^{1/2}$$

$$+ \int_0^t C_{\mathbf{k},n} e^{-\lambda(t-s)} \| F_n(s) - \langle F_n(s)\rangle \|_{|\mathbf{k}|, \beta_{n,|\mathbf{k}|}} ds (1+|x|^{\alpha_{\mathbf{k},n}}).$$

Using the induction assumption, we have for $t \geq 0$, $i \in \mathbb{N}^*$ and $k \in \mathbb{N}^d$

$$|\partial_\mathbf{k}(v_n(t) - \int_{\mathbb{R}^d} \phi(y)\mu_n(y)\rho(y)dy)| \leq \sum_{m=1}^n c_m P_{0,n-m}(t) e^{-\lambda t} \| \phi - \langle \phi \rangle \|_{4n, l_{0,4n}}$$

$$+ \sum_{\ell=0}^{n-1} \int_0^t C_{\mathbf{k},\ell} \tilde{P}_{\mathbf{k},\ell}(s) e^{-\lambda(t-s)} e^{-\lambda s} ds \| \phi - \langle \phi \rangle \|_{|\mathbf{k}|+4n, \ell_{|\mathbf{k}|,n}} (1+|x|^{\alpha_{\mathbf{k},n}}).$$

The conclusion follows. ♦

The following Proposition ends the proof of our main result Theorem 2.9.

**Proposition 5.4** *Let $N$ and $\ell_N$ be fixed. Let $\delta_0 = \frac{1}{\alpha}$. Let $X_p$ be the discrete process defined by the implicit Euler scheme (2.10) or the implicit split-step scheme (2.9). Let $0 \leq \delta < \delta_0$ and $\phi \in \mathcal{C}^\infty_{pol}(\mathbb{R}^d) \cap \mathcal{C}^{6N+2}_{\ell_N}(\mathbb{R}^d)$. Then there exist constants $C_N$ and $p_N$ such that for all $p \in \mathbb{N}$,*

$$\| \mathbb{E}\phi(X_p) - v^{(N)}(t_p,.) \|_{0,p_N} \leq C_N \| \phi - \langle \phi \rangle \|_{6N+2,\ell_N} \delta^N, \tag{5.4}$$

*where $t_p = p\delta$.*
*Moreover, we have for $0 \leq \delta < \delta_0$ and for all function $\phi \in \mathcal{C}^\infty_{pol}(\mathbb{R}^d) \cap \mathcal{C}^{6N+2}_{\ell_N}(\mathbb{R}^d)$ that there exist constants $C_N$ and $p_N$ and a positive polynomial function $P_N$ satisfying the following : For all $p \in \mathbb{N}$,*

$$\| \mathbb{E}\phi(X_p) - \int_{\mathbb{R}^d} \phi(x)\mu^N(x)\rho(x)dx \|_{0,p_N} \leq \| \phi - \langle \phi \rangle \|_{6N+2,\ell_N} \Big( e^{-\lambda t_p} P_N(t_p) + C_N \delta^N \Big),$$

*where $t_p = p\delta$.*

**Proof.** Let $N$ and $\ell_N$ be fixed. Let $\delta_0 = \frac{1}{\alpha}$. Let $X_p$ be the discrete process defined by the implicit Euler scheme (2.10) or the implicit split-step scheme (2.9). Let $0 \leq \delta < \delta_0$ and $\phi \in \mathcal{C}^\infty_{pol}(\mathbb{R}^d) \cap \mathcal{C}^{6N+2}_{\ell_N}(\mathbb{R}^d)$. For all $p$, where $t_j = j\delta$ for $j \leq p$, we have for $x \in \mathbb{R}^d$ such that $X_0 = x$

$$\mathbb{E}\phi(X_p) - v^{(N)}(t_p, x) = \mathbb{E}v^{(N)}(0, X_p) - v^{(N)}(t_p, x)$$
$$= \mathbb{E}\sum_{j=0}^{p-1} \mathbb{E}^{X_{p-j-1}}\Big(v^{(N)}(t_j, X_{p-j}) - v^{(N)}(t_{j+1}, X_{p-j-1})\Big).$$

Here we have used the notation $\mathbb{E}^{X_{p-j-1}}$ for the conditional expectation with respect to the filtration generated by $X_{p-j-1}$. We have:

$$\mathbb{E}^{X_{p-j-1}}\Big(v^{(N)}(t_j, X_{p-j}) - v^{(N)}(t_{j+1}, X_{p-j-1})\Big)$$
$$= \mathbb{E}^{X_{p-j-1}}\Big(v^{(N)}(t_j, X_1(X_{p-j-1})) - v^{(N)}(t_{j+1}, X_{p-j-1})\Big),$$

where $X_1(x)$ is the first step of the scheme (2.10) or of the scheme (2.9) when the initial condition is $x$. Using Theorem 4.1 with $t = t_j$, Proposition 2.5 and (5.2), we deduce that there exist integers $p_N$ and $k_N$ such that

$$\| \mathbb{E}\phi(X_p) - v^{(N)}(t_p, .) \|_{0, p_N} \leq \delta^{N+1} C_N \sum_{j=0}^{p-1} \sup_{\substack{s \in ]0,\delta[ \\ n=0,\ldots,N}} |v_n(t_{j+1}, .)|_{2N+2, k_N}$$

$$\leq \delta^{N+1} C_N \sum_{j=0}^{p-1} e^{-\lambda t_j} P_N(t_j) \| \phi - \langle\phi\rangle \|_{6N+2, \ell_N}$$

$$\leq \delta^{N+1} C_N \| \phi - \langle\phi\rangle \|_{6N+2, \ell_N} \sum_{j=0}^{p-1} e^{-\tilde{\lambda} t_j},$$

for some constant $C_N$. We have used:

$$|v_n(t_{j+1}, .)|_{2N+2, k_N} = |v_n(t_{j+1}, .) - \int_{\mathbb{R}^d} \phi(x)\mu_n(x)\rho(x)dx|_{2N+2, k_N}.$$

We conclude by using the fact that for a fixed constant $\hat\gamma > 0$, we have

$$\sum_{j=0}^{p-1} e^{-\hat\gamma j \delta} \leq \frac{1}{1 - e^{-\hat\gamma \delta}} \leq \frac{C}{\delta},$$

where the constant C depends on $\tilde\gamma$ and $\delta_0$. This shows (5.4). The second estimate is a consequence of (5.2). ♦

## A  Appendix : Proof of Proposition 2.7

We warn the reader that the constants may vary from line to line during the proofs, and that in order to use lighter notations we usually forget to mention dependence on the parameters. We use the generic notation $C$ for such constants.

The aim of this appendix is to show the following result: Let $V \in C^\infty(\mathbb{R}^d)$ such that $V$ verifies the assumptions **B**.

**Proposition A.1** Let $\phi \in \mathcal{C}_{pol}^{\infty}(\mathbb{R}^d)$ such that $\int_{\mathbb{R}^d} \phi(x)\rho(x)dx = 0$. Let $u$ be the solution of

$$\frac{d}{dt}u(t,x) = L(x)u(t,x), \ x \in \mathbb{R}^d, \ t > 0, \quad u(0,x) = \phi(x), \ x \in \mathbb{R}^d, \tag{A.1}$$

where $L$ is defined for all $x \in \mathbb{R}^d$ and $\phi \in \mathcal{C}_{pol}^{\infty}(\mathbb{R}^d)$ by

$$L(x)\phi(x) = \frac{1}{2}\sum_{i=1}^{d} \partial_{ii}\phi(x) - \sum_{i=1}^{d} \partial_i V(x)\partial_i\phi(x).$$

For any integer $m$, there exist integers $\ell_m$ and $s > \ell_m$ and strictly positive real numbers $C$ and $\lambda$ such that $\phi \in \mathcal{C}_{\ell_m}^m(\mathbb{R}^d)$ and for all $t > 0$, $\mathbf{k} \in \mathbb{N}^d$ such that $|\mathbf{k}| = m$ and $x \in \mathbb{R}^d$

$$|\partial_{\mathbf{k}} u(t,x)| \leq C(1+|x|^s)\exp(-\lambda t) \parallel \phi \parallel_{m,\ell_m}. \tag{A.2}$$

The proof of Proposition A.1 proceeds as follow

1. We first show that $u \in \mathcal{C}_{pol}^{\infty}(\mathbb{R}^d)$ and the result (A.2) for all $t \leq 1$.
2. We show the point-wise estimate for $u$.
3. Then, using the Bismuth-Elworthy formulas, we show the result (A.2) for $t = 1$.
4. Finally, using the last two items, we show the result (A.2) for $t \geq 1$.

We recall that we use the following notation: Let $x \in \mathbb{R}^d$ be fixed, $X_x(t)$ is the solution of (2.1) such that $X(0) = x$.

## A.1 The polynomial growth of $u$ and its derivatives

**Lemma A.2** Let $\phi \in \mathcal{C}_{pol}^{\infty}(\mathbb{R}^d)$. The function $u$, defined by (2.16), and all its derivatives have polynomial growth: For all $p$, there exist some constants $s_p, \ell_p \in \mathbb{N}$, $\gamma_p$ and $C$ such that $\phi \in \mathcal{C}_{2\ell_p}^p(\mathbb{R}^d)$ and for all $x \in \mathbb{R}^d$, $\mathbf{k} \in \mathbb{N}^d$ such that $|\mathbf{k}| = p$ and $t > 0$,

$$|\partial_{\mathbf{k}} u(t,x)| \leq C \exp(\gamma_p t)(1+|x|^{2s_p}) \parallel \phi \parallel_{p,2\ell_p}. \tag{A.3}$$

**Proof.** In all the proof, $C$ is an ever changing constant.
Let us show the result (A.3) for $p = 0$. Let us assume that $\phi \in \mathcal{C}_{2\ell_0}^0(\mathbb{R}^d)$. Using Proposition 2.2, we have for all $x \in \mathbb{R}^d$ and $t > 0$

$$u(t,x) = \mathbb{E}(\phi(X_x(t)))$$
$$|u(t,x)| \leq \parallel \phi \parallel_{0,2\ell_0} \left(\mathbb{E}|X_x(t)|^{2\ell_0} + 1\right)$$
$$\leq C \parallel \phi \parallel_{0,2\ell_0} (|x|^{2\ell_0} + 1).$$

Let us now show the result (A.3) for $p = 1$. We have for all $x \in \mathbb{R}^d$, $h \in \mathbb{R}^d$ and $t > 0$

$$Du(t,x) \cdot h = \mathbb{E}\big(D\phi(X_x(t)) \cdot \eta_x^h(t)\big), \tag{A.4}$$

where $\eta_x^h(t) \in \mathbb{R}^d$ is a process defined for $x \in \mathbb{R}^d$ and $h \in \mathbb{R}^d$ by

$$\eta_x^h(t) = DX_x(t) \cdot h \quad \text{for } t > 0$$

and $\eta_x^h(0) = h$. Moreover, we have for all $t > 0$, $x \in \mathbb{R}^d$ and $h \in \mathbb{R}^d$

$$\frac{d}{dt}\eta_x^h(t) = -D^2 V(X_x(t)) \cdot \eta_x^h(t).$$

By definition of $\phi$, we have that there exists $\ell_1 \in \mathbb{N}$ such that $\phi \in \mathcal{C}_{\ell_1}^1(\mathbb{R}^d)$. Then, using Proposition 2.2 and (A.4), we have for $x \in \mathbb{R}^d$, $h \in \mathbb{R}^d$ and $t > 0$

$$\begin{aligned}|Du(t,x) \cdot h| &\leq \| \phi \|_{1,\ell_1} \left(\mathbb{E}(|X_x(t)|^{\ell_1}|\eta_x^h(t)|) + \mathbb{E}(|\eta_x^h(t)|)\right)\\ &\leq \| \phi \|_{1,\ell_1} \left((\mathbb{E}|X_x(t)|^{2\ell_1}\mathbb{E}|\eta_x^h(t)|^2)^{1/2} + (\mathbb{E}|\eta_x^h(t)|^2)^{1/2}\right)\\ &\leq C \| \phi \|_{1,\ell_1} (|x|^{\ell_1} + 1)(\mathbb{E}|\eta_x^h(t)|^2)^{1/2}.\end{aligned}$$

Moreover, using semi-convexity assumption B-2, we have for all $x \in \mathbb{R}^d$, $h \in \mathbb{R}^d$ and $t > 0$

$$\frac{d}{dt}|\eta_x^h(t)|^2 = -2D^2 V(X_x(t)) \cdot (\eta_x^h(t), \eta_x^h(t)) \leq 2\alpha |\eta_x^h(t)|^2,$$

where $\alpha$ is the constant of semi-convexity of $V$. Using Gronwall's lemma, we obtain for all $x \in \mathbb{R}^d$, $h \in \mathbb{R}^d$ and $t > 0$

$$|\eta_x^h(t)|^2 \leq e^{2\alpha t}|h|^2 \tag{A.5}$$

and

$$|Du(t,x) \cdot h| \leq C \| \phi \|_{1,\ell_1} \exp(2\alpha t)(|x|^{\ell_1} + 1)|h|.$$

Then, we have for all $x \in \mathbb{R}^d$, $t > 0$ and $i \in \{1, ..., d\}$

$$|\partial_i u(t,x)| \leq C \| \phi \|_{1,\ell_1} \exp(2\alpha t)(|x|^{\ell_1} + 1).$$

Let us show Lemma A.2 for $p = 2$. We have for $x \in \mathbb{R}^d$, $h \in \mathbb{R}^d$ and $t > 0$

$$D^2 u(t,x) \cdot (h,h) = \mathbb{E}\Big(D^2\phi(X_x(t)) \cdot (\eta_x^h(t), \eta_x^h(t)) + D\phi(X_x(t)) \cdot \xi_x^h(t)\Big),$$

where $\xi_x^h(t) \in \mathbb{R}^d$ is a process defined for $x \in \mathbb{R}^d$, $h \in \mathbb{R}^d$ by

$$\xi_x^h(t) = D^2 X_x(t) \cdot (h,h) \quad \text{for } t > 0,$$

and $\xi_x^h(0) = 0$. Moreover, we have for $x \in \mathbb{R}^d$, $h \in \mathbb{R}^d$ and $t > 0$

$$\frac{d}{dt}\xi_x^h(t) = -D^3 V(X_x(t)) \cdot (\eta_x^h(t), \eta_x^h(t)) - D^2 V(X_x(t)) \cdot \xi_x^h(t).$$

Using assumption B-4 on the polynomial growth of $V$, we have that there exists $p_1 \in \mathbb{N}$ such that $V \in \mathcal{C}_{2p_1}^3(\mathbb{R}^d)$. Using semi-convexity assumption B-2, we have for $x \in \mathbb{R}^d$, $h \in \mathbb{R}^d$ and $t > 0$

$$\begin{aligned}\frac{d}{dt}|\xi_x^h(t)|^2 &= -2\big(D^3 V(X_x(t)) \cdot (\eta_x^h(t), \eta_x^h(t), \xi_x^h(t)) - D^2 V(X_x(t)) \cdot (\xi_x^h(t), \xi_x^h(t))\big)\\ &\leq 2 \| V \|_{3, 2p_1} (|X_x(t)|^{2p_1} + 1)|\eta_x^h(t)|^2 |\xi_x^h(t)| + 2\alpha |\xi_x^h(t)|^2\\ &\leq C_1 \Big[(|X_x(t)|^{2p_1} + 1)^2 |\eta_x^h(t)|^4 + |\xi_x^h(t)|^2\Big].\end{aligned}$$

Using Gronwall's lemma and Proposition 2.2, we obtain that there exists a constant $\tilde{\gamma}$ depending of $C_1$ and of the constant of semi-convexity $\alpha$ such that for all $x \in \mathbb{R}^d$, $h \in \mathbb{R}^d$ and $t > 0$

$$\begin{aligned}\mathbb{E}|\xi_x^h(t)|^2 &\leq C\exp(4\alpha t)|h|^4(|x|^{4p_1} + 1) + (|x|^{4p_1} + 1)|h|^4 \int_0^t C\exp(4\alpha s)\exp((2\alpha + C_1)s)ds\\ &\leq C\exp(\tilde{\gamma}t)(|x|^{4p_1} + 1)|h|^4,\end{aligned}$$

where $C$ is an ever changing constant. Moreover, we have that there exists $\ell_2 \in \mathbb{N}$ such that $\phi \in \mathcal{C}^2_{\ell_2}(\mathbb{R}^d)$. Finally, we have that there exist some constants $k_2$ and $\gamma_2$ such that for all $x \in \mathbb{R}^d$, $h \in \mathbb{R}^d$ and $t > 0$

$$|D^2 u(t,x) \cdot (h,h)| \leq \| \phi \|_{2,\ell_2} \mathbb{E}((|X_x(t)|^{\ell_2} + 1)|\eta_x^h(t)|^2) + \| \phi \|_{2,\ell_2} \left( \mathbb{E}(|X_x(t)|^{2\ell_2} + 1)\mathbb{E}|\xi_x^h(t)|^2 \right)^{1/2}$$

$$\leq C \| \phi \|_{2,\ell_2} \exp(\gamma_2 t)(1 + |x|^{2k_2})|h|^2.$$

This show the result (A.3) for $p = 2$.

To show the result for higher derivatives, we use an induction, the Faà di Bruno's formula and the same methods used for the first derivatives. Moreover, we can prove that for all $k \geq 0$, there exist some constants $q_k$, $\alpha_k$ and $C$ such that for all $x \in \mathbb{R}^d$, $h \in \mathbb{R}^d$ and $t > 0$

$$\mathbb{E}|D^{k+1} X_x(t) \cdot (h,...,h)|^2 \leq C \exp(\alpha_k t)(1 + |x|^{2q_k})^2 |h|^{2(k+1)}. \tag{A.6}$$

♦

**Remark A.3** *The result (A.2) for $t \leq 1$ is a corollary of this Lemma.*

## A.2 Estimate of $u$

**Proposition A.4** *Let $\phi \in \mathcal{C}^\infty_{pol}(\mathbb{R}^d)$ be fixed such that $\int_{\mathbb{R}^d} \phi(x)\rho(x)dx = 0$ and $u$ be the solution of (A.1). Let us assume that $\phi \in \mathcal{C}^0_{2\ell_0}(\mathbb{R}^d)$. There exist a positive real numbers $C$ and $\lambda$ such that for all $x \in \mathbb{R}^d$ and $t \geq 0$*

$$|u(t,x)| \leq C \| \phi \|_{0,2\ell_0} \exp(-\lambda t)(1 + |x|^{2\ell_0}).$$

A proof of this result can be found in [11]. Our equation is dissipative and has noise in all the direction, then Proposotion A.4 is a corollary of Theorem 4.4 of [11].

## A.3 Estimate of the derivatives of $u$

We can now show an estimate of the derivatives of $u$ at the time $t = 1$:

**Lemma A.5** *Let $\phi \in \mathcal{C}^\infty_{pol}(\mathbb{R}^d)$ and $k \in \mathbb{N}^*$ be fixed and $u$ be the solution of (A.1). Let us assume that $\phi \in \mathcal{C}^0_{2\ell_0}(\mathbb{R}^d)$. There exist constants $C$ and $m_k \geq 2\ell_0$ such that we have for all $x \in \mathbb{R}^d$ and $\mathbf{j} \in \mathbb{N}^d$ such that $|\mathbf{j}| = k$*

$$|\partial_{\mathbf{j}} u(1,x)| \leq C \| \phi \|_{0,2\ell_0} (1 + |x|^{m_k}).$$

The Lemma A.5 is a corollary of the following lemma.

**Lemma A.6** *Let $\phi \in \mathcal{C}^\infty_{pol}(\mathbb{R}^d) \cap \mathcal{C}^0_{2\ell_0}(\mathbb{R}^d)$ be fixed and $u$ be the solution of (A.1). For $k \in \mathbb{N}^*$, there exist constants $C$ and $m_k \geq 2\ell_0$ such that we have for $0 < t \leq 1$, $x \in \mathbb{R}^d$ and $\mathbf{j} \in \mathbb{N}^d$ such that $|\mathbf{j}| = k$*

$$|\partial_{\mathbf{j}} u(t,x)| \leq C \| \phi \|_{0,2\ell_0} t^{-k/2}(1 + |x|^{m_k}). \tag{A.7}$$

**Proof.** We only prove the result for the two first derivatives, as the result for the higher order follows from analogous arguments and an induction.

Let us show the result (A.7) for $k = 1$. We have the Bismuth-Elworthy formula (see [5]): for $x \in \mathbb{R}^d$, $h \in \mathbb{R}^d$ and $0 < t \leq 1$

$$Du(t,x) \cdot h = \frac{1}{t} \mathbb{E} \Big( u(0, X_x(t)) \int_0^t \langle \eta_x^h(s), dW(s) \rangle \Big),$$

where $\eta_x^h(t) = DX_x(t) \cdot h \in \mathbb{R}^d$ for $t > 0$ and $\eta_x^h(0) = h$.

Using Cauchy-Schwarz inequality, Burkholder-Davis-Gundy inequality and (A.5) which bounds $\eta_x^h(t)$, we have for $t \leq 1$, $x \in \mathbb{R}^d$ and $h \in \mathbb{R}^d$

$$\begin{aligned}|Du(t,x) \cdot h| &\leq t^{-1}(\mathbb{E}(|u(0, X_x(t))|^2)^{1/2}(\mathbb{E}(\int_0^t |\eta_x^h(s)|^2 ds))^{1/2}) \\ &\leq t^{-1/2}(\mathbb{E}(|\phi(X_x(t))|^2)^{1/2} C|h|).\end{aligned}$$

Using $\phi \in \mathcal{C}_{2\ell_0}^0(\mathbb{R}^d)$ and Proposition 2.2 on the moment of the solution of (2.1), we get (A.7) for k=1.

Let us show it for $k = 2$. We have the Bismuth-Elworthy at the second order (see [5]): for $x \in \mathbb{R}^d$, $h \in \mathbb{R}^d$ and $t > 0$

$$\begin{aligned}D^2 u(t,x) \cdot (h,h) =& \frac{2}{t}\Big[\mathbb{E}\Big(Du(t/2, X_x(t/2)) \cdot (DX_x(t/2) \cdot h) \int_0^{t/2} \langle DX_x(t/2) \cdot h, dW(s)\rangle\Big) \\ &+ \mathbb{E}\Big(\int_0^{t/2} Du(t-s, X_x(s)) \cdot (D^2 X_x(s) \cdot (h,h)) ds\Big)\Big]. \end{aligned} \quad (A.8)$$

Using Cauchy-Schwarz inequality, Burkholder-Davis-Gundy inequality, (A.6) and (A.5) which bound $DX_x(t) \cdot h$ and $D^2 X_x(t) \cdot (h,h)$, we get for $0 < t \leq 1$, $x \in \mathbb{R}^d$ and $h \in \mathbb{R}^d$

$$\begin{aligned}&|D^2 u(t,x) \cdot (h,h)| \\ &\leq \frac{2}{t}\Big(\mathbb{E}|Du(t/2, X_x(t/2))|^4\Big)^{1/4} \Big(\mathbb{E}|DX_x(t/2) \cdot h|^4\Big)^{1/4} \Big(\mathbb{E}(\int_0^{t/2} |DX_x(t/2) \cdot h|^2 ds)\Big)^{1/2} \\ &\quad + \frac{2}{t}\int_0^{t/2} \Big(\mathbb{E}|Du(t-s, X_x(s))|^2\Big)^{1/2} \Big(\mathbb{E}|D^2 X_x(s) \cdot (h,h)|^2\Big)^{1/2} ds \\ &\leq \frac{C}{t}(1+|x|^q)|h|^2\Big[t^{1/2}\Big(\mathbb{E}|Du(t/2, X_x(t/2))|^4\Big)^{1/4} + \int_0^{t/2} \Big(\mathbb{E}|Du(t-s, X_x(s))|^2\Big)^{1/2} ds\Big].\end{aligned}$$

To conclude, we use the result (A.7) for $k = 1$ and the Proposition 2.2 on the moment of the solution of the equation (2.1).

Using the proof of (A.8) done in [5], we show the following formula at the order 3 for $x \in \mathbb{R}^d$, $h \in \mathbb{R}^d$ and $t > 0$

$$\begin{aligned}&D^3 u(t,x) \cdot (h,h,h) \\ &= \frac{2}{t}\Big[\mathbb{E}\Big(D^2 u(t/2, X_x(t/2)) \cdot (DX_x(t/2) \cdot h, DX_x(t/2) \cdot h) \int_0^{t/2} \langle DX_x(s) \cdot h, dW(s)\rangle\Big) \\ &\quad + \mathbb{E}\Big(Du(t/2, X_x(t/2)) \cdot \big(D^2 X_x(t/2) \cdot (h,h)\big) \int_0^{t/2} \langle DX_x(s) \cdot h, dW(s)\rangle\Big) \\ &\quad + 2\mathbb{E}\Big(\int_0^{t/2} D^2 u(t-s, X_x(s)) \cdot \big(D^2 X_x(s) \cdot (h,h), DX_x(s) \cdot h\big) ds\Big) \\ &\quad + \mathbb{E}\Big(\int_0^{t/2} Du(t-s, X_x(s)) \cdot \big(D^3 X_x(s) \cdot (h,h,h)\big) ds\Big)\Big].\end{aligned}$$

We use the same arguments as previously to prove (A.7) for $k = 3$. The proof for highest order use the same ideas. ♦

We can know show the result (A.2) for $t \geq 1$.

**Lemma A.7** *Let $\phi \in \mathcal{C}_{pol}^\infty(\mathbb{R}^d)$ be fixed such that $\phi \in \mathcal{C}_{2\ell_0}^0(\mathbb{R}^d)$ and $u$ be the solution of (A.1). For all $k \in \mathbb{N}^*$, there exist constant $C_k$ and $m_k$ such that for all $t \geq 1$, $\mathbf{j} = (j_1, ..., j_d) \in \mathbb{N}^d$ such that $|\mathbf{j}| = k$ and $x \in \mathbb{R}^d$,*

$$|\partial_{\mathbf{j}} u(t,x)| \leq C_k \exp(-\lambda t)(1+|x|^{m_k}) \| \phi \|_{0,2\ell_0},$$

*where $\lambda$ is defined in Proposition A.4.*

**Remark A.8** *The constant $m_k$ depends of the polynomial growth of all the derivatives of $\phi$ of order less that $k$ and the polynomial growth of all the derivatives of $V$ of order less that $k+1$.*

**Proof.** Let $x \in \mathbb{R}^d$ and $t \geq 1$ be fixed. Let us assume that $\phi \in \mathcal{C}_{2\ell_0}^0(\mathbb{R}^d)$. We have for all $t \geq 1$,

$$u(t,x) = \mathbb{E}(u(t-1, X_x(1))) = \mathbb{E}\bigl(P_{t-1}(X_x(1))\bigr).$$

Let $v$ defined for all $s > 0$ by $v(s,x) = \mathbb{E}P_{t-1}(X_x(s))$, then $u(t,x) = v(1,x)$. Using Proposition A.4, we have for all $y \in \mathbb{R}^d$

$$|P_{t-1}(y)| \leq C(1+|y|^{2\ell_0}) \exp(-\lambda t) \| \phi \|_{0,2\ell_0}.$$

For $k \in \mathbb{N}^*$, using Lemma A.6 on $v$ and $t \geq 1$, we get there exist constants $m_k \geq 2\ell_0$ and $C$ such that for any $\mathbf{j} \in \mathbb{N}^d$ such that $|\mathbf{j}| = k$

$$\begin{aligned}|\partial_{\mathbf{j}} u(t,x)| &\leq C \| P_{t-1} \|_{0,2\ell_0} (1+|x|^{m_k}) \\ &\leq C(1+|x|^{m_k}) \exp(-\lambda t) \| \phi \|_{0,2\ell_0}.\end{aligned}$$

♦

# References

[1] ABDULLE, A., COHEN D., VILMART, G. & ZYGALAKIS, K.(2012) High order weak methods for stochastic differential equations based on modified equations. *SIAM J.Sci. Comput.* 34 : 1800-1823.

[2] CERRAI,S. (2001) Second order PDE's in finite and infinite dimension: A probabilistic approach, volume 1762 of *Lecture Notes in Mathematics*. Springer-Verlag, Berlin.

[3] DEBUSSCHE, A. & FAOU, E. (2009) Modified energy for split-step methods applied to the linear Schrödinger equation. *SIAM J. Numer. Anal.*, 47(5):3705–3719.

[4] DEBUSSCHE, A. & FAOU, E. (2012) Weak backward error analysis for SDEs. *SIAM J. Numer. Anal.*

[5] ELWORTHY, K. D. & LI, X. M. (1994) Formulae for the derivatives of heat semigroups. *J. Funct. Anal.*, 125(1):252–286.

[6] FAOU, E. (2012) Geometric numerical integration and Schrödinger equations. *Zurich Lectures in Advanced Mathematics. European Mathematical Society (EMS)*, Zürich.

[7] FAOU, E. & GRÉBERT, B. (2011) Hamiltonian interpolation of splitting approximations for nonlinear PDEs. *Found. Comput. Math.*, 11(4):381–415.

[8] FRIEDMAN, A. (1975) Stochastic differential equations and applications. Vol. 1. Academic Press [Harcourt Brace Jovanovich Publishers], New York. *Probability and Mathematical Statistics*, Vol. 28.

[9] HAIRER, E. & LUBICH, CH. (1997) The life-span of backward error analysis for numerical integrators. *Numer. Math.*, 76(4):441–462.

[10] HAIRER, E., LUBICH, CH. & WANNER, G. (2010) Geometric numerical integration. Structure-preserving algorithms for ordinary differential equations, volume 31 of *Springer Series in Computational Mathematics*. Springer, Heidelberg. Reprint of the second (2006) edition.

[11] HIGHAM, D. J., MATTINGLY, J. C. & STUART A. M. (2002) Ergodicity for SDEs and approximations: locally Lipschitz vector fields and degenerate noise. *Stochastic Process. Appl.*, 101(2):185–232.

[12] KLOEDEN, P. E. & PLATEN, E. (1992) Numerical solution of stochastic differential equations, volume 23 of *Applications of Mathematics (New York)*. Springer-Verlag, Berlin.

[13] LEIMKUHLER, B. & REICH, S. (2004) Simulating Hamiltonian dynamics. *Cambridge Monographs on Applied and Computational Mathematics*, 14. Cambridge University Press, Cambridge.

[14] LELIÈVRE, T., ROUSSET, M. & STOLTZ, G. (2010) Free energy computations. *Imperial College Press*, London, A mathematical perspective.

[15] LIONS, J. L. (1969) Quelques méthodes de résolution des problèmes aux limites non linéaires. *Dunod*.

[16] MATTINGLY, J. C., STUART, A. M. & TRETYAKOV, M. V. (2010) Convergence of numerical time-averaging and stationary measures via Poisson equations. *SIAM J. Numer. Anal.*, 48(2):552–577.

[17] MILSTEIN, G. N. & TRETYAKOV, M. V. (2004) Stochastic numerics for mathematical physics. *Scientific Computation*. Springer-Verlag, Berlin.

[18] MOSER, J. K. (1968) Lectures on Hamiltonian systems. *Mem. Amer. Math. Soc.* No. 81. Amer. Math. Soc., Providence, R.I.

[19] PARDOUX, E & VERETENNIKOV A. YU. (2001) On the Poisson equation and diffusion approximation. I. *Ann. Probab.*, 29(3):1061–1085.

[20] REICH, S., (1999) Backward error analysis for numerical integrators, *SIAM J. Numer. Anal.* 36, 1549–1570.

[21] ROYER, G. (1999) Une initiation aux inégalités de Sobolev logarithmiques, volume 5 of *Cours Spécialisés [Specialized Courses]*. Société Mathématique de France, Paris.

[22] SHARDLOW, T. (2006) Modified equations for stochastic differential equations. *BIT*, 46(1):111–125.

[23] TALAY, D. (1990) Second order discretization schemes of stochastic differential systems for the computation of the invariant low. *Stochastics and Stochastic Reports*, 29(1):13–36.

[24] TALAY, D. (1996) Probabilistic numerical methods for partial differential equations: elements of analysis. In *Probabilistic models for nonlinear partial differential equations (Montecatini Terme, 1995)*, volume 1627 of *Lecture Notes in Math.*, pages 148–196. Springer, Berlin.

[25] TALAY, D. (2002) Stochastic Hamiltonian systems: exponential convergence to the invariant measure, and discretization by the implicit Euler scheme. *Markov Process. Related Fields*, 8(2):163–198.

[26] TALAY, D. & TUBARO, L. (1990) Expansion of the global error for numerical schemes solving stochastic differential equations. *Stochastic Anal. Appl.*, 8(4):483–509.